\documentclass[a4paper,12pt]{article}
\usepackage{amsmath, amsfonts, amssymb}
\usepackage[all]{xy}
\makeatletter
\newbox\bk@bxb
\newbox\bk@bxa
\newif\if@bkcont
\newcount\bk@lcnt

\def\breakboxskip{2pt}
\def\breakboxparindent{1.8em}

\def\breakbox{\vskip\breakboxskip\relax
\setbox\bk@bxb\vbox\bgroup
\advance\linewidth -2\fboxrule
\hsize\linewidth\@parboxrestore
\parindent\breakboxparindent\relax}

\def\bk@split{%
\@tempdimb\ht\bk@bxb 
\advance\@tempdimb\dp\bk@bxb
\setbox\bk@bxa\vsplit\bk@bxb to\z@ 
\setbox\bk@bxa\vbox{\unvbox\bk@bxa}
\setbox\@tempboxa\vbox{\copy\bk@bxa\copy\bk@bxb}
\advance\@tempdimb-\ht\@tempboxa
\advance\@tempdimb-\dp\@tempboxa}

\def\bk@addfsepht{%
\setbox\bk@bxa\vbox{\vskip\fboxsep\box\bk@bxa}}

\def\bk@addskipht{%
\setbox\bk@bxa\vbox{\vskip\@tempdimb\box\bk@bxa}}

\def\bk@addfsepdp{%
\@tempdima\dp\bk@bxa
\advance\@tempdima\fboxsep
\dp\bk@bxa\@tempdima}

\def\bk@addskipdp{%
\@tempdima\dp\bk@bxa
\advance\@tempdima\@tempdimb
\dp\bk@bxa\@tempdima}

\def\bk@line{%
\hbox to \linewidth{%
\hskip-2\fboxsep\vrule \@width\fboxrule\hskip.5\fboxsep\vrule \@width\fboxrule\hskip1.5\fboxsep
\box\bk@bxa\hfil
}}%

\def\endbreakbox{\egroup
\ifhmode\par\fi{\noindent\bk@lcnt\@ne
\@bkconttrue\baselineskip\z@\lineskiplimit\z@
\lineskip\z@\vfuzz\maxdimen
\bk@split\bk@addfsepht\bk@addskipdp
\ifvoid\bk@bxb 
\def\bk@fstln{\bk@addfsepdp
\hskip-\parindent\vbox{\llap{\raisebox{-2ex}{\rule{1.5\fboxsep}{\fboxrule}\hskip.5\fboxsep}}\bk@line\llap{\rule{1.5\fboxsep}{\fboxrule}\hskip.5\fboxsep}}}

\else 
\def\bk@fstln{\vbox{\llap{\raisebox{-2ex}{\rule{1.5\fboxsep}{\fboxrule}\hskip.5\fboxsep}}\bk@line}\hfil%
\advance\bk@lcnt\@ne
\loop
\bk@split\bk@addskipdp\leavevmode
\ifvoid\bk@bxb 
\@bkcontfalse\bk@addfsepdp
\vtop{\bk@line\llap{\rule{2\fboxsep}{\fboxrule}}}%

\else 
\bk@line
\fi
\hfil\advance\bk@lcnt\@ne
\if@bkcont\repeat}%
\fi
\leavevmode\bk@fstln\par}\vskip\breakboxskip\relax}

\newcommand{\fracb}[2]{\frac{\raisebox{-.7ex}{$\scriptstyle #1$}}{\raisebox{.7ex}{$\scriptstyle #2$}}}

\def\smp{\smallskip\par}
\def\un{{\bf 1}}
\def\zero{\{0\}}
\def\pf{\noindent{\bf Proof~:}\ }

\def\findemo{~\leaders\hbox to 1em{\hss\  \hss}\hfill~\raisebox{.5ex}{\framebox[1ex]{}}\smp}

\def\mpn{\medskip\par\noindent}
\def\smpn{\smallskip\par\noindent}
\def\normal{\mathop{\trianglelefteq}}

\def\smp{\smallskip\par}
\def\smpn{\smallskip\par\noindent}

\def\mpoint{\;\;.}
\def\mvirg{\;\;,}

\def\Res{{\rm Res}}
\def\Ind{{\rm Ind}}
\def\Inf{{\rm Inf}}
\def\Def{{\rm Def}}
\def\Iso{{\rm Iso}}

\def\Indinf{{\rm Indinf}}
\def\Defres{{\rm Defres}}

\def\Hom{{\rm Hom}}

\def\Ext{{\rm Ext}}
\def\Inf{{\rm Inf}}

\def\Out{{\rm Out}}
\def\Aut{{\rm Aut}}

\def\Ker{{\rm Ker}}

\def\Id{{\rm Id}}

\def\op{^{op}}

\def\dsp{\displaystyle}
\def\Z{\mathbb{Z}}
\def\N{\mathbb{N}}

\newcommand{\dirsum}[1]{\mathop{\oplus}_{#1}\limits}
\newcommand{\romain}[1]{\uppercase\expandafter{\romannumeral #1}}

\newcommand{\flh}[2]{\mathop{\hbox to 12mm{\rightarrowfill}}_{\displaystyle #2}^{\displaystyle #1}\limits}
\newcommand{\sflh}[2]{\mathop{\hbox to 12mm{\rightarrowfill}}_{\scriptstyle #2}^{\scriptstyle #1}\limits}

\newcommand{\gMod}[1]{#1{\hbox{-}\mathsf{Mod}}}

\newcommand{\sou}[1]{\underline{#1}}
\newcommand{\sur}[1]{\,\overline{\! #1}}
\newcommand{\sumb}[2]{\mathop{\sum}_{{\scriptstyle #1}\atop {\scriptstyle #2}}}

\newcommand{\sumc}[3]{\sum_{{\scriptstyle #1}\atop {{\scriptstyle #2}\atop {\scriptstyle #3}}}}

\def\op{^{op}}

\newcommand{\carre}[8]{\begin{array}{ccc}
#1&\mathop{\hbox to 12mm{\rightarrowfill}}^{\displaystyle{#2}}\limits&#3\\
\llap{$\displaystyle{#4}$}\left\downarrow\vbox to 6mm{}\right. & & \left\downarrow\vbox to 6mm{}\right.\rlap{$\displaystyle{#5}$}\\
#6&\mathop{\hbox to 12mm{\rightarrowfill}}_{\displaystyle #7}\limits&#8\\
\end{array}}

\newcommand{\carrem}[8]{\begin{array}{ccc}
#1&\mathop{\hbox to 12mm{\rightarrowfill}}^{\displaystyle #2}\limits&#3\\
\llap{$\displaystyle #4$}\left\uparrow\vbox to 6mm{}\right. & & \left\uparrow\vbox to 6mm{}\right.\rlap{$\displaystyle #5$}\\
#6&\mathop{\hbox to 12mm{\rightarrowfill}}_{\displaystyle #7}\limits&#8\\
\end{array}}

\newcommand{\limproj}[1]{\lim_{\displaystyle\longleftarrow\atop \scriptstyle{#1}}\limits}

\newenvironment{enonce}[1]{\pagebreak[2]\refstepcounter{subsection}\refstepcounter{prop}\smpn{{\bf \thesection.\arabic{prop}.\ #1:}}\begin{it} }{\end{it}\smp}
\newenvironment{enonce*}[1]{\pagebreak[2]\smpn{#1:}\begin{it} }{\end{it}\smp}
\newcommand{\result}[1]{\begin{enonce}{#1}}
\def\fresult{\end{enonce}}
\newcommand{\npar}{\smallskip\par\noindent\pagebreak[2]\refstepcounter{subsection}\refstepcounter{prop}{\bf \thesection.\arabic{prop}.\ \ }}

\newenvironment{mth}[1]{\begin{breakbox}\begin{enonce}{#1}}{\end{enonce}\end{breakbox}}
\newenvironment{mth*}[1]{\begin{breakbox}\begin{enonce*}{#1}}{\end{enonce*}\end{breakbox}}
\newenvironment{rem}[1]{\refstepcounter{subsection}\refstepcounter{prop} \mpn{{\bf \thesection.\arabic{prop}.}\ \bf#1:}}{\smp}

\def\dom{\backslash}
\makeatletter
\renewenvironment{enumerate}{\ifnum \@enumdepth >3 \@toodeep\else
      \advance\@enumdepth \@ne
      \edef\@enumctr{enum\romannumeral\the\@enumdepth}\list
      {\csname label\@enumctr\endcsname}{\setlength{\topsep}{1ex}\setlength{\itemsep}{0pt}\usecounter
        {\@enumctr}\def\makelabel##1{\hss\llap{##1}}}\fi}{\endlist}
\renewenvironment{itemize}{\ifnum \@itemdepth >3 \@toodeep\else \advance\@itemdepth \@ne
\edef\@itemitem{labelitem\romannumeral\the\@itemdepth}%
\list{\csname\@itemitem\endcsname}{\setlength{\topsep}{1ex}\setlength{\itemsep}{0pt}\def\makelabel##1{\hss\llap{##1}}}\fi}
{\endlist}
\makeatother
\makeatletter
\def\@sect#1#2#3#4#5#6[#7]#8{\ifnum #2>\c@secnumdepth
    \let\@svsec\@empty\else
    \refstepcounter{#1}\edef\@svsec{\csname the#1\endcsname .\hskip .5em}\fi
    \@tempskipa #5\relax
     \ifdim \@tempskipa>\z@
       \begingroup #6\relax
         \@hangfrom{\hskip #3\relax\@svsec}{\interlinepenalty \@M #8\par}%
       \endgroup
      \csname #1mark\endcsname{#7}\addcontentsline
        {toc}{#1}{\ifnum #2>\c@secnumdepth \else
                     \protect\numberline{\csname the#1\endcsname}\fi
                   #7}\else
       \def\@svsechd{#6\hskip #3\relax  
                  \@svsec #8\csname #1mark\endcsname
                     {#7}\addcontentsline
                          {toc}{#1}{\ifnum #2>\c@secnumdepth \else
                            \protect\numberline{\csname the#1\endcsname}\fi
                      #7}}\fi
    \@xsect{#5}}
\def\section{\@startsection {section}{1}{\z@}{-3.5ex plus-1ex minus
    -.2ex}{2.3ex plus.2ex}{\center\reset@font\large\bf}}  

\makeatother
\renewenvironment{equation}{\refstepcounter{subsection}\refstepcounter{prop}$$}{\leqno{\bf (\theprop)}$$}

\def\mar[#1]{\ar@{-}[#1]|-{\object@{<}}}
\def\marb[#1]{\ar@{-}[#1]|{\object+{  }}}

\def\CC{\mathcal{C}}
\def\CM{\mathcal{M}}
\def\CE{\mathcal{E}}
\def\CF{\mathcal{F}}
\def\CG{\mathcal{G}}
\def\CI{\mathcal{I}}
\def\CJ{\mathcal{J}}
\def\CO{\mathcal{O}}
\def\CR{\mathcal{R}}
\def\CS{\mathcal{S}}

\def\CY{\mathcal{Y}}
\def\Cf{R}
\def\At{\mathcal{A}t}
\def\El{\CE{l}}
\def\endpf{\findemo}
\newcommand{\enr}[1]{[L]}
\begin{document}
\centerline{\Large\bf Idempotents of double Burnside algebras,}\vspace{2ex}
\centerline{\Large\bf $L$-enriched bisets,}\vspace{2ex}
\centerline{\Large\bf and decomposition of $p$-biset functors}\vspace{1cm}\par
\centerline{\bf Serge Bouc}\vspace{.1cm}\par
\begin{center}
\begin{minipage}{11cm}
{\footnotesize {\bf Abstract:} Let $\Cf$ be a (unital) commutative ring, and $G$ be a finite group with order invertible in $\Cf$. We introduce new idempotents $\epsilon_{T,S}^G$ in the double Burnside algebra $\Cf B(G,G)$ of $G$ over~$\Cf$,  
indexed by conjugacy classes of minimal sections $(T,S)$ of $G$ (i.e. sections such that $S\leq\Phi(T)$). These idempotents are orthogonal, and their sum is equal to the identity. It follows that for any biset functor $F$ over~$\Cf$, the evaluation $F(G)$ splits as a direct sum of specific $\Cf$-modules indexed by minimal sections of $G$, up to conjugation.\par
The restriction of these constructions to the biset category of $p$-groups, where $p$ is a prime number invertible in $\Cf$, leads to a decomposition of the category of $p$-biset functors over $\Cf$ as a direct product of categories $\mathcal{F}_L$ indexed by {\em atoric} $p$-groups $L$ up to isomorphism. \par
We next introduce the notions of {\em $L$-enriched biset} and {\em $L$-enriched biset functor} for an arbitrary finite group $L$, and show that for an atoric $p$-group $L$, the category $\CF_L$ is equivalent to the category of $L$-enriched biset functors defined over elementary abelian $p$-groups.\par
Finally, the notion of {\em vertex} of an indecomposable $p$-biset functor is introduced (when $p\in \Cf^\times$), and when $\Cf$ is a field of characteristic different from $p$, the objects of the category $\CF_L$ are characterized in terms of vertices of their composition factors.\par
}\vspace{1ex}
\vspace{1ex}\par
\end{minipage}
\end{center}
\section{Introduction}
Let $\Cf$ denote throughout a commutative ring (with identity element). For a finite group $G$,  we consider the double Burnside algebra $\Cf B(G,G)$ of a $G$ over $\Cf$. In the case where the order of $G$ is invertible in $\Cf$, we introduce idempotents $\epsilon_{T,S}^G$ in $\Cf B(G,G)$, indexed by the set $\mathcal{M}(G)$ of minimal sections of $G$, i.e. the set of pairs $(T,S)$ of subgroups of $G$ with $S\normal T$ and $S\leq \Phi(T)$, where $\Phi(T)$ is the Frattini subgroup of $G$ (such sections have been considered in Section~5 of~\cite{vanishing}). The idempotent $\epsilon_{T,S}^G$ only depends of the conjugacy class of $(T,S$) in $G$. Moreover, the idempotents $\epsilon_{T,S}^G$, where $(T,S)$ runs through a set $[\mathcal{M}(G)]$ of representatives of orbits of $G$ acting on $\mathcal{M}(G)$ by conjugation, are orthogonal, and their sum is equal to the identity element of $\Cf B(G,G)$. \par
The idempotents $\epsilon_{G,\un}^G$ plays a special role in our construction, and it is denoted by $\varphi_\un^G$. In particular, when $F$ is a biset functor over $\Cf$ (and the order of $G$ is invertible in $\Cf$), we set $\delta_\Phi F(G)=\varphi_\un^GF(G)$. We show that $\delta_\Phi F(G)$ consists of those elements $u\in F(G)$ such that $\Res_H^Gu=0$ whenever $H$ is a proper subgroup of $G$, and $\Def_{G/N}^Gu=0$ whenever $N$ is a non-trivial normal subgroup of $G$ contained in $\Phi(G)$. This yields moreover a decomposition 
$$F(G)\cong\big(\dirsum{(T,S)\in\mathcal{M}(G)}\delta_\Phi F(T/S)\big)^G\cong\dirsum{(T,S)\in[\mathcal{M}(G)]}\delta_\Phi F(T/S)^{N_G(T,S)/T}\mpoint$$
Restricting these constructions to the biset category $\Cf\CC_p$ of $p$-groups with coefficients in $\Cf$, where $p$ is a prime invertible in~$\Cf$, we get orthogonal idempotents $b_L$ in the {\em center} of $\Cf\CC_p$, indexed by {\em atoric} $p$-groups, i.e. finite $p$-groups which cannot be split as a direct product $C_p\times Q$, for some $p$-group $Q$. We show next that every finite $p$-group $P$ admits a unique largest atoric quotient~$P^@$, well defined up to isomorphism, and that there exists an elementary abelian $p$-subgroup $E$ of $P$ (non unique in general) such that $P\cong E\times P^@$. For a given atoric $p$-group $L$, we introduce a category $\Cf\CC_p^{\sharp L}$, defined as a quotient of the subcategory of $\Cf\CC_p$ consisting of $p$-groups $P$ such that $P^@\cong L$. This leads to a decomposition of the category $\CF_{p,\Cf}$ of $p$-biset functors over $\Cf$ as a direct product
$$\CF_{p,\Cf}\cong\prod_{L\in[\At_p]}\mathsf{Fun}_\Cf\big(\Cf\CC_p^{\sharp L},\gMod{\Cf}\big)$$
 of categories of representations of $\Cf\CC_p^{\sharp L}$ over $\Cf$, where $L$ runs through a set $[\At_p]$ of isomorphism classes of atoric $p$-groups. Similar questions on idempotents in double Burnside algebras and decomposition of biset functors categories have been considered by L. Barker~(\cite{barker-blocks-Mackey}), R. Boltje and S. Danz (\cite{boltje-danz-ghost}, \cite{boltje-danz-ghost-zero}), R. Boltje and B. K\"ulshammer (\cite{boltje-kulshammer}), and P. Webb (\cite{webb-stratification-II}).\par
In particular, via the above decomposition, to any indecomposable $p$-biset functor $F$ is associated a unique atoric $p$-group, called the {\em vertex} of $F$. We show that this vertex is isomorphic to $Q^@$, for any $p$-group $Q$ such that $F(Q)\neq \zero$ but $F$ vanishes on any proper subquotient of $Q$.\par
Going back to arbitrary finite groups, we next introduce the notions of {\em $L$-enriched biset} and {\em $L$-enriched biset functor}, and show that when $L$ is an atoric $p$-group, the abelian category $\mathsf{Fun}_\Cf\big(\Cf\CC_p^{\sharp L},\gMod{\Cf}\big)$ is equivalent to the category of $L$-enriched biset functors from elementary abelian $p$-groups to $\Cf$-modules.\par
The paper is organized as follows: Section~2 is a review of definitions and basic results on Burnside rings and biset functors. Section~3 is concerned with the algebra~$\mathcal{E}(G)$ obtained by ``cutting'' the double Burnside algebra $\Cf B(G,G)$ of a finite group $G$ by the idempotent $\widetilde{e_G^G}$ corresponding to the ``top'' idempotent $e_G^G$ of the Burnside algebra $\Cf B(G)$. Orthogonal idempotents $\varphi_N^G$ of $\mathcal{E}(G)$ are introduced, indexed by normal subgroups $N$ of $G$ contained in $\Phi(G)$. It is shown moreover that if $G$ is nilpotent, then $\varphi_\un^G$ is central in $\mathcal{E}(G)$. In Section~4, the idempotents $\epsilon_{T,S}^G$ of $\Cf B(G,G)$ are introduced, leading in Section~5 to the corresponding direct sum decomposition of the evaluation at $G$ of any biset functor over $\Cf$.
In Section~6, atoric $p$-groups are introduced, and their main properties are stated. In Section~7, the biset category of $p$-groups over $\Cf$ is considered, leading to a splitting of the category $\CF_{p,\Cf}$ of $p$-biset functors over $\Cf$ as a direct product of abelian categories $\CF_L=\mathsf{Fun}_\Cf\big(\Cf\CC_p^{\sharp L},\gMod{\Cf}\big)$ indexed by atoric $p$-groups $L$ up to isomorphism. In Section 8, for an arbitrary finite group $L$, the notions of $L$-enriched biset and $L$-enriched biset functor are introduced, and it is shown that when $L$ is an atoric $p$-group, the category $\CF_L$ is equivalent to the category of $L$-enriched biset functors on elementary abelian $p$-groups. Finally, in Section~9, for a given atoric $p$-group $L$, and when $p$ is invertible in $\Cf$, the structure of the category $\CF_L$ is considered, and the notion of vertex of an indecomposable $p$-biset functor over $\Cf$ is introduced. In particular, when $\Cf$ is a field of characteristic different from $p$, it is shown that the objects of $\CF_L$ are those $p$-biset functors all composition factors of which have vertex $L$.
\tableofcontents
\section{Review of Burnside rings and biset functors}
\npar Let $G$ be a finite group, let $s_G$ denote the set of subgroups of $G$, let $\sur{s_G}$ denote the set of conjugacy classes of subgroups of $G$, and let $[s_G]$ denote a set of representatives of $\sur{s_G}$. \par
Let $B(G)$ denote the Burnside ring of $G$, i.e. the Grothendieck ring of the category of finite $G$-sets. It is a commutative ring, with an identity element, equal to the class of a $G$-set of cardinality 1. The additive group $B(G)$ is a free abelian group on the set $\{[G/H]\mid H\in[s_G]\}$ of isomorphism classes of transitive $G$-sets.\par
\npar $\bullet$ When $G$ and $H$ are finite groups, and $L$ is a subgroup of $G\times H$, set
\begin{eqnarray*}
p_1(L)&=&\{g\in G\mid\exists h\in H,\;(g,h)\in L\}\mvirg\\
p_2(L)&=&\{h\in H\mid\exists g\in G,\;(g,h)\in L\}\mvirg\\
k_1(L)&=&\{g\in G\mid(g,1)\in L\}\mvirg\\
k_2(L)&=&\{h\in H\mid(1,h)\in L\}\mpoint\\
\end{eqnarray*}
Recall that $k_i(L)\normal p_i(L)$, for $i\in\{1,2\}$, that $\big(k_1(L)\times k_2(L)\big)\normal L$, and that there are canonical isomorphisms
$$p_1(L)/k_1(L)\cong L/\big(k_1(L)\times k_2(L)\big)\cong p_2(L)/k_2(L)\mpoint$$
Set moreover $q(L)=L/\big(k_1(L)\times k_2(L)\big)$.\mpn
$\bullet$ When $Z$ is a subgroup of $G$, set
$$\Delta(Z)=\{(z,z)\mid z\in Z\}\leq (G\times G)\mpoint$$
When $N$ is a normal subgroup of $G$, set
$$\Delta_N(G)=\{(a,b)\in G\times G\mid ab^{-1}\in N\}\mpoint$$
It is a subgroup of $G\times G$.\mpn
$\bullet$ When $G$, $H$, and $K$ are groups, when $L\leq (G\times H)$ and $M\leq (H\times K)$, set
$$L*M=\{(g,k)\in (G\times K)\mid\exists h\in H,\;(g,h)\in L\;\hbox{and}\;(h,k)\in K\}\mpoint$$
It is a subgroup of $(G\times K)$.
\npar When $G$ and $H$ are finite groups, a $(G,H)$-{\em biset} $U$ is a set endowed with a left action of $G$ and a right action of $H$ which commute. In other words $U$ is a $G\times H\op$-set, where $H\op$ is the opposite group of $H$. The {\em opposite biset} $U\op$ is the $(H,G)$-biset equal to $U$ as a set, with actions defined for $h\in H$, $u\in U$ and $g\in G$ by $h\cdot u\cdot g\; \hbox{(in $U\op$)}=g^{-1}uh^{-1}\; \hbox{(in $U$)}$. \par
The Burnside group $B(G,H)$ is the Grothendieck group of the category of finite $(G,H)$-bisets. It is a free abelian group on the set of isomorphism classes $[(G\times H)/L]$, for $L\in [s_{G\times H}]$, where the $(G,H)$-biset structure on $(G\times H)/L$ is given by
$$\forall a,g\in G,\,\forall b,h\in H,\;a\cdot(g,h)L\cdot b=(ag,b^{-1}h)L\mpoint$$
When $G$, $H$, and $K$ are finite groups, there is a unique bilinear product
$$\times_H:B(G,H)\times B(H,K)\to B(G,K)$$
induced by the usual product $(U,V)\mapsto U\times_HV=(U\times V)/H$ of bisets, where the right action of $H$ on $U\times V$ is defined for $u\in U$, $v\in V$ and $h\in H$ by $(u,v)\cdot h=(uh,h^{-1}v)$. This product will also be denoted as a composition $(\alpha,\beta)\mapsto \alpha\circ\beta$ or as a product $(\alpha,\beta)\mapsto \alpha\beta$.\par
\pagebreak[3]
This leads to the following definitions:
\pagebreak[3]
\begin{mth}{Definition} The {\em biset category of finite groups} $\mathcal{C}$ is defined as follows:
\begin{itemize}
\item The objects of $\mathcal{C}$ are the finite groups.
\item When $G$ and $H$ are finite groups, 
$$\Hom_\mathcal{C}(G,H)=B(H,G)\mpoint$$
\item When $G$, $H$, and $K$ are finite groups, the composition
$$\circ : \Hom_\mathcal{C}(H,K)\times \Hom_\mathcal{C}(G,H)\to\Hom_\mathcal{C}(G,K)$$
is the product 
$$\times_H: B(K,H)\times B(H,G)\to B(K,G)\mpoint$$
\item The identity morphism of the group $G$ is the class of the set $G$, viewed as a $(G,G)$-biset by left and right multiplication.
\end{itemize}
A {\em biset functor} is an additive functor from $\mathcal{C}$ to the category of abelian groups. 
\end{mth}
When $\Cf$ is a commutative (unital) ring, the category $\Cf\mathcal{C}$ is defined similarly by extending coefficients to $\Cf$, i.e. by setting
$$\Hom_{\Cf\mathcal{C}}(G,H)=\Cf\otimes_\Z B(H,G)\mvirg$$
which will be simply denoted by $\Cf B(H,G)$. A {\em biset functor over $\Cf$} is an $\Cf$-linear functor from $\Cf\mathcal{C}$ to the category $\gMod{\Cf}$ of $\Cf$-modules. The category of biset functors over $\Cf$ (where morphisms are natural transformations of functors) is denoted by $\CF_\Cf$.\par
The correspondence sending a $(G,H)$-biset $U$ to its opposite $U\op$ extends to an isomorphism of $\Cf$-modules $\Cf B(G,H)\to \Cf B(H,G)$. These isomorphisms give an equivalence of $\Cf$-linear categories from $\Cf\mathcal{C}$ to its opposite category, which is the identity on objects.
\npar Let $G$ and $H$ be finite groups, and $F$ be a biset functor (with values in $\gMod{\Cf}$). For any finite $(H,G)$-biset $U$, the isomorphism class $[U]$ of $U$ belongs to $B(H,G)$, and it yields an $\Cf$-linear map $F([U]): F(G)\to F(H)$, simply denoted by $F(U)$, or even $f\in F(G)\mapsto U(f)\in F(H)$. In particular:
\begin{itemize}
\item When $H$ is a subgroup of $G$, denote by $\Ind_H^G$ the set $G$, viewed as a $(G,H)$-biset for left and right multiplication, and by $\Res_H^G$ the same set, viewed as an $(H,G)$-biset. This gives a map $\Ind_H^G:F(H)\to F(G)$, called induction, and a map $\Res_H^G:F(G)\to F(H)$, called restriction.
\item When $N$ is a normal subgroup of $G$, let $\Inf_{G/N}^G$ denote the set $G/N$, viewed as a $(G,G/N)$-biset for the left action of $G$, and right action of $G/N$ by multiplication. Also let $\Def_{G/N}^G$ denote the set $G/N$, viewed as a $(G/N,G)$-biset. This gives a map $\Inf_{G/N}^G:F(G/N)\to F(G)$, called inflation, and a map $\Def_{G/N}^G:F(G)\to F(G/N)$, called deflation.
\item Finally, when $f:G\to G'$ is a group isomorphism, denote by $\Iso(f)$ the set $G'$, viewed as a $(G',G)$-biset for left multiplication in $G'$, and right action of $G$ given by multiplication by the image under $f$. This gives a map $\Iso(f):F(G)\to F(G')$, called transport by isomorphism.
\end{itemize}
When $G$ and $H$ are finite groups, any $(G,H)$-biset is a disjoint union of transitive ones. It follows that any element of $B(G,H)$ is a linear combination of morphisms of the form $[(G\times H)/L]$, where $L\in s_{G\times H}$. Moreover, any such morphism factors as
\begin{equation}\label{factorize}[(G\times H)/L]=\Ind_{p_1(L)}^G\circ\Inf_{p_1(L)/k_1(L)}^{p_1(L)}\circ\Iso(f)\circ\Def_{p_2(L)/k_2(L)}^{p_2(L)}\circ\Res_{p_2(L)}^H\mvirg
\end{equation}
where $f:p_2(L)/k_2(L)\to p_1(L)/k_1(L)$ is the canonical group isomorphism.\par
In particular, for $N\normal G$,
\begin{equation}\label{infdef}
[(G\times G)/\Delta_N(G)]=\Inf_{G/N}^G\circ\Def_{G/N}^G\mpoint
\end{equation}
For finite groups $G, H, K$, for $L\leq (G\times H)$ and $M\leq (H\times K)$, one has that
\begin{equation}\label{Mackey bisets}[(G\times H)/L]\times_H[(H\times K)/M]=\sum_{h\in p_2(L)\dom H/p_1(M)}[(G\times K)/(L*{^{(h,1)}M})
\end{equation}
in $B(G,K)$.
\npar When $G$ is a finite group, the group $B(G,G)$ is the ring of endomorphisms of $G$ in the category $\mathcal{C}$. This ring is called the double Burnside ring of $G$. It is a non-commutative ring (if $G$ is non trivial), with identity element equal to the class of the set~$G$, viewed as a $(G,G)$-biset for left and right multiplication.\par
There is a unitary ring homomorphism $\alpha\mapsto \widetilde{\alpha}$ from $B(G)$ to $B(G,G)$, induced by the functor $X\mapsto \widetilde{X}$ from $G$-sets to $(G,G)$-bisets, where $\widetilde{X}=G\times X$, with $(G,G)$-biset structure given by
$$\forall a,b,g\in G,\,\forall x\in X,\;a(g,x)b=(agb,ax)\mpoint$$
This construction has in particular the following properties (\cite{bisetfunctors}, Corollary 2.5.12):
\pagebreak[3]
\begin{mth}{Lemma} \label{tilde compatibility}Let $G$ be a finite group.\begin{enumerate}
\item If $H$ is a subgroup of $G$, and $X$ is a finite $G$-set, then there is an isomorphism of $(G,H)$-bisets
$$\widetilde{X}\times_G\Ind_H^G\cong \Ind_H^G\times_H\widetilde{\Res_H^GX}\mvirg$$
and an isomorphism of $(H,G)$-bisets
$$\Res_H^G\times_G\widetilde{X}\cong \widetilde{\Res_H^GX}\times_H\Res_H^G\mpoint$$
\item If $H$ is a subgroup of $G$, and $Y$ is a finite $H$-set, then there is an isomorphism of $(G,G)$-bisets
$$\Ind_H^G\times_H\widetilde{Y}\times_H\Res_H^G\cong\widetilde{\Ind_H^GY}\mpoint$$
\item If $N$ is a normal subgroup of $G$, and $X$ is a finite $G/N$-set, then there is an isomorphism of $(G/N,G)$-bisets
$$\widetilde{X}\times_{G/N}\Def_{G/N}^G\cong \Def_{G/N}^G\times_G\widetilde{\Inf_{G/N}^GX}\mpoint$$
\item If $N$ is a normal subgroup of $G$, and $X$ is a finite $G$-set, then there is an isomorphism of $(G/N,G/N)$-bisets
$$\Def_{G/N}^G\times_GX\times_G\Inf_{G/N}^G\cong \widetilde{\Def_{G/N}^GX}\mpoint$$
\end{enumerate}
\end{mth}
\npar Let $\Cf B(G)$ denote the $\Cf$-algebra $\Cf\otimes_\Z B(G)$. Assume moreover that the order of $G$ is invertible in $\Cf$. Then for $H\leq G$, let $e_H^G\in\Cf B(G)$ be defined by
\begin{equation}\label{ehg}
e_H^G=\frac{1}{|N_G(H)|}\sum_{K\leq H}|K|\mu(K,H)\;[G/K]\mvirg
\end{equation}
where $\mu$ is the M\"obius function of the poset of subgroups of $G$. The elements $e_H^G$, for $H\in[s_G]$, are orthogonal idempotents of $\Cf B(G)$, and their sum is equal to the identity element of $\Cf B(G)$. It follows that the elements $\widetilde{e_H^G}$, for $H\in[s_G]$, are orthogonal idempotents of the $\Cf$-algebra $\Cf B(G,G)=\Cf\otimes_\Z B(G,G)$, and the sum of these idempotents is equal to the identity element of $\Cf B(G,G)$. The idempotents $\widetilde{e_G^G}$ play a special role, due to the following lemma:
\pagebreak[3]
\begin{mth}{Lemma} \label{general}Let $\Cf$ be a commutative ring, and $G$ be a finite group with order invertible in $\Cf$.
\begin{enumerate}
\item Let $H$ be a proper subgroup of $G$. Then
$$\Res_H^G\circ \widetilde{e_G^G}=0\;\;\hbox{and}\;\;\widetilde{e_G^G}\circ\Ind_H^G=0\mpoint$$
\item When $N\normal G$, let $m_{G,N}\in\Cf$ be defined by
$$m_{G,N}=\frac{1}{|G|}\sumb{X\in s_G}{XN=G}|X|\mu(X,G)\mpoint$$
Then
$$\Def_{G/N}^G\circ \widetilde{e_G^G}\circ\Inf_{G/N}^G=m_{G,N}\widetilde{e_{G/N}^{G/N}}\mpoint$$
\item Let $N\normal G$, and suppose that $N$ is contained in the Frattini subgroup $\Phi(G)$ of $G$. Then
$$\widetilde{e_{G/N}^{G/N}}\circ\Def_{G/N}^G=\Def_{G/N}^G\circ\widetilde{e_G^G}\;\;\hbox{and}\;\;\Inf_{G/N}^G\circ\widetilde{e_{G/N}^{G/N}}=\widetilde{e_G^G}\circ\Inf_{G/N}^G\mpoint$$
\end{enumerate}
\end{mth}
\pf Assertion 1 follows from Lemma~\ref{tilde compatibility} and Assertion 1 of Theorem~5.2.4. of \cite{bisetfunctors}. \par
Assertion 2 follows from  Lemma~\ref{tilde compatibility} and Assertion 4 of Theorem~5.2.4. of \cite{bisetfunctors}.\par
Finally the first part of Assertion 3 follows from Lemma~\ref{tilde compatibility} and Assertion 3 of Theorem~5.2.4. of \cite{bisetfunctors}: indeed $\Inf_{G/N}^Ge_{G/N}^{G/N}$ is equal to the sum of the different idempotents $e_X^G$ of $\Cf B(G)$ indexed by subgroups $X$ such that $XN=G$. If $N\leq \Phi(G)$, then $XN=G$ implies $X\Phi(G)=G$, hence $X=G$. The second part of Assertion~3 follows by taking opposite bisets, since $\widetilde{e_G^G}$ and $\widetilde{e_{G/N}^{G/N}}$ are equal to their opposite bisets, and since $(\Def_{G/N}^G)\op\cong \Inf_{G/N}^G$.\findemo
\begin{rem}{Remark} \label{mgn Frattini} For the same reason, if $N\leq \Phi(G)$, then $m_{G,N}=1$.
\end{rem}
\begin{rem}{Remark} \label{saturant} It follows from Assertion 1 and Remark~\ref{factorize} that if $G$ and $H$ are finite groups and if $L\leq (G\times H)$, then $\widetilde{e_G^G}[(G\times H)/L]=0\rule{0ex}{3ex}$ if $p_1(L)\neq G$, and $[(G\times H)/L]\widetilde{e_H^H}=0$ if $p_2(L)\neq H$.
\end{rem}
\section{Idempotents in $\mathcal{E}(G)$}
\begin{mth}{Notation} When $G$ is a finite group with order invertible in $\Cf$, denote by $\mathcal{E}(G)$ the $\Cf$-algebra defined by
$$\mathcal{E}(G)=\widetilde{e_G^G}\Cf B(G,G)\widetilde{e_G^G}\mpoint$$
Set
$$\Sigma(G,G)=\{M\in s_{G\times G}\mid p_1(L)=p_2(L)=G\}\mvirg$$
and for $L\in s_{G\times G}$, set
$$Y_L=\widetilde{e_G^G}\,[(G\times G)/L]\,\widetilde{e_G^G}\in\mathcal{E}(G)\mpoint$$
\end{mth}
The $\Cf$-algebra $\mathcal{E}(G)$ has been considered in~\cite{doublact}, Section~9, in the case $\Cf$ is a field of characteristic 0. The extension of the results proved there to the case where $\Cf$ is a commutative ring in which the order of $G$ is invertible is straightforward. In particular:
\begin{mth}{Proposition} \label{product of YL}Let $G$ be a finite group with order invertible in $\Cf$.
\begin{enumerate}
\item If $L\in s_{G\times G}-\Sigma(G,G)$, then $Y_L=0$.
\item The elements $Y_L$, for $L$ in a set of representatives of $(G\times G)$-conjugacy classes on $\Sigma(G,G)$, form a $\Cf$-basis of $\mathcal{E}(G)$.
\item For $L,M\in\Sigma(G,G)$
$$Y_LY_M=\frac{m_{G,k_2(L)\cap k_1(M)}}{|G|}\sumc{Z\leq G}{\rule{0ex}{1.5ex}Zk_2(L)=Zk_1(M)=G}{\rule{0ex}{1.5ex}Z\geq k_2(L)\cap k_1(M)}|Z|\mu(Z,G)\;Y_{L*\Delta(Z)*M}$$
in $\mathcal{E}(G)$.
\end{enumerate}
\end{mth}
\begin{mth}{Corollary} \label{sub Frattini}Let $L, M\in\Sigma(G,G)$. If one of the groups $k_2(L)$ or $k_1(M)$ is contained in $\Phi(G)$, then
$$Y_LY_M=Y_{L*M}\mpoint$$
\end{mth}
\pf Indeed if $k_2(L)\leq\Phi(G)$, then $Zk_2(L)=G$ implies $Z\Phi(G)=G$, hence $Z=G$. Similarly, if $k_1(M)\leq\Phi(G)$, then $Zk_1(M)=G$ implies $Z=G$. In each case, Proposition~\ref{product of YL} then gives
$$Y_LY_M=m_{G,k_2(L)\cap k_1(M)}Y_{L*M}\mvirg$$
and moreover $m_{G,k_2(L)\cap k_1(M)}=1$ since $k_2(L)\cap k_1(M)\leq\Phi(G)$, by Remark~\ref{mgn Frattini}.\endpf
\begin{mth}{Notation} For a normal subgroup $N$ of $G$ such that $N\leq \Phi(G)$, set
$$\varphi_N^G=\sumb{M\normal G}{N\leq M\leq\Phi(G)}\mu_{\normal G}(N,M)Y_{\Delta_M(G)}\mvirg$$
where $\mu_{\normal G}$ is the M\"obius function of the poset of normal subgroups of $G$.
\end{mth}
\vspace{3ex}
\pagebreak[3]
\begin{mth}{Proposition} \label{phi1}Let $N\normal G$ with $N\leq\Phi(G)$. Then
$$\varphi_N^G=\Inf_{G/N}^G\varphi_\un^{G/N}\Def_{G/N}^G\mpoint$$
\end{mth}
\pf Indeed if $N\leq M\normal G$, then $\mu_{\normal G}(N,M)=\mu_{\normal G/N}(\un,M/N)$. Since moreover $N\leq \Phi(G)$, setting $\sur{G}=G/N$ and $\sur{M}=M/N$, we have by Lemma~\ref{general}
\begin{eqnarray*}
\Inf_{G/N}^GY_{\Delta_{\sur{G}}(\sur{M})}\Def_{G/N}^G&=&\Inf_{G/N}^G\circ\widetilde{e_{\sur{G}}^{\sur{G}}}\big((\sur{G}\times\sur{G})/\Delta_{\sur{G}}(\sur{M})\big)\widetilde{e_{\sur{G}}^{\sur{G}}}\circ\Def_{G/N}^G\\
&=&\widetilde{e_G^G}\circ\Inf_{G/N}^G\big((\sur{G}\times\sur{G})/\Delta_{\sur{G}}(\sur{M})\big)\Def_{G/N}^G\circ\widetilde{e_G^G}\\
&=&\widetilde{e_G^G}\big((G\times G)/\Delta_M(G)\big)\widetilde{e_G^G}\\
&=&Y_{\Delta_M(G)}\mvirg
\end{eqnarray*}
since $\Inf_{G/N}^G\big((\sur{G}\times\sur{G}))/\Delta_{\sur{G}}(\sur{M})\big)\Def_{G/N}^G=(G\times G)/\Delta_M(G)$.\endpf
\pagebreak[3]
\begin{mth}{Proposition} \label{phinG}\begin{enumerate}
\item Let $N\normal G$, with $N\leq\Phi(G)$. Then
\begin{eqnarray*}
\varphi_N^G&=&\widetilde{e_G^G}\times_G\Big(\sumb{M\normal G}{N\leq M\leq\Phi(G)}\mu_{\normal G}(N,M)[(G\times G)/\Delta_M(G)]\Big)\\
&=&\Big(\sumb{M\normal G}{N\leq M\leq\Phi(G)}\mu_{\normal G}(N,M)[(G\times G)/\Delta_M(G)]\Big)\times_G\widetilde{e_G^G}\mpoint
\end{eqnarray*}
\item In particular
$$\varphi_\un^G=\frac{1}{|G|}\sumb{X\leq G, M\normal G}{\rule{0ex}{1.6ex}M\leq\Phi(G)\leq X\leq G} |X|\mu(X,G)\mu_{\normal G}(\un,M)\,\Indinf_{X/M}^G\circ\Defres_{X/M}^G\mpoint$$
\end{enumerate}
\end{mth}
\pf For Assertion 1, by definition
$$\varphi_N^G=\hspace{-3ex}\sumb{M\normal G}{N\leq M\leq \Phi(G)}\hspace{-3ex}\mu_{\normal G}(N,M)\widetilde{e_G^G}[(G\times G)/\Delta_M(G)]\times_G\sum_{X\leq G}\fracb{|X|}{|G|}\mu(X,G)[(G\times G)/\Delta(X)].$$
Moreover $[(G\times G)/\Delta_M(G)]\times_G[(G\times G)/\Delta(X)]=[(G\times G)/\big(\Delta_M(G)*\Delta(X)\big)]$, by~(\ref{Mackey bisets}), and $\Delta_M(G)*\Delta(X)=\{(xm,x)\mid x\in X, m\in M\}$. The first projection of this group is equal to $XM$, hence it is equal to $G$ if and only if $X=G$, since $M\leq \Phi(G)$. The first equality of Assertion~1 follows, by Remark~\ref{saturant}. The second one follows by taking opposite bisets, since $\widetilde{e_G^G}$ and $[(G\times G)/\Delta_M(G)]$ are equal to their opposite.\par
Assertion 2 follows in the special case where $N=1$, expanding $\widetilde{e_G^G}$ as
$$\widetilde{e_G^G}=\frac{1}{|G|}\sum_{X\leq G}|X|\mu(X,G)[(G\times G)/\Delta(X)]\mvirg$$
observing that $\mu(X,G)=0$ unless $X\geq \Phi(G)$, and that if $X\geq \Phi(G)\geq M$, then
$$[(G\times G)/\Delta(X)]\circ[(G\times G)/\Delta_M(G)]=[(G\times G)/\Delta_M(X)]\mvirg$$
which is equal to $\Indinf_{X/M}^G\circ\Defres_{X/M}^G$.  \endpf
\pagebreak[3]
\begin{mth}{Corollary}\label{essentiel frattini} \begin{enumerate} 
\item Let $H<G$. Then $\Res_H^G\varphi_N^G=0$ and $\varphi_N^G\Ind_H^G=0$.
\item Let $M\normal G$. If $M\cap \Phi(G)\nleq N$, then $\Def_{G/M}^G\varphi_N^G=0$ and $\varphi_N^G\Inf_{G/M}^G=0$.
\end{enumerate}
\end{mth}
\pf The first part of Assertion 1 follows from Lemma~\ref{general}, since 
$$\Res_H^G\varphi_N^G=\Res_H^G\widetilde{e_G^G}\varphi_N^G=0\mpoint$$
The second part follows by taking opposite bisets.\par
For Assertion 2, let $P=M\cap \Phi(G)$. Since $\Def_{G/M}^G=\Def_{G/M}^{G/P}\circ\Def_{G/P}^G$, it suffices to consider the case $M=P$, i.e. the case where $M\leq \Phi(G)$. Then, since $[(G\times G)/\Delta_M(G)]=\Inf_{G/M}^G\Def_{G/M}^G$ for any $M\normal G$, by~\ref{infdef}, and since $\Def_{G/M}^G\Inf_{G/Q}^G=\Inf_{G/MQ}^{G/M}\Def_{G/MQ}^{G/Q}$ for any $M,Q\normal G$,
\begin{eqnarray*}
\Def_{G/M}^G\varphi_N^G&=&\Def_{G/M}^G\sumb{Q\normal G}{\rule{0ex}{1.5ex}N\leq Q\leq\Phi(G)}\mu_{\normal G}(N,Q)\Inf_{G/Q}^G\Def_{G/Q}^G\widetilde{e_G^G}\\
&=&\sumb{Q\normal G}{\rule{0ex}{1.5ex}N\leq Q\leq\Phi(G)}\mu_{\normal G}(N,Q)\Inf_{G/MQ}^G\Def_{G/MQ}^G\widetilde{e_G^G}\\
&=&\sumb{P\normal G}{\rule{0ex}{1.5ex}N\leq P\leq \Phi(G)}\big(\sumc{Q\normal G}{\rule{0ex}{1.5ex}N\leq Q\leq \Phi(G)}{\rule{0ex}{1.5ex}QM=P}\mu_{\normal G}(N,Q)\big)\Inf_{G/P}^G\Def_{G/P}^G\widetilde{e_G^G}\mpoint
\end{eqnarray*}
Now for a given $P\normal G$ with $P\subseteq \Phi(G)$, the sum $\sumc{Q\normal G}{\rule{0ex}{1.5ex}N\leq Q\leq \Phi(G)}{\rule{0ex}{1.5ex}QM=P}\limits\mu_{\normal G}(N,Q)$ is equal to zero unless $NM=N$, that is $M\leq N$, by classical properties of the M\"obius function (\cite{stanley} Corollary 3.9.3). This proves the first part of Assertion~2, and the second part follows by taking opposite bisets.\endpf
\pagebreak[3]
\begin{mth}{Theorem} \label{idempotents E}Let $G$ be a finite group with order invertible in~$\Cf$.\begin{enumerate}
\item The elements $\varphi_N^G$, for $N\normal G$ with $N\leq\Phi(G)$, form a set of orthogonal idempotents in the algebra $\mathcal{E}(G)$, and their sum is equal to the identity element $\widetilde{e_G^G}$ of $\mathcal{E}(G)$.
\item Let $N\normal G$ with $N\leq \Phi(G)$, and let $H$ be a finite group. \begin{enumerate}
\item If $L\leq (G\times H)$, then $\varphi_N^G\times_G[(G\times H)/L]=0$ unless $p_1(L)=G$ and $k_1(L)\cap\Phi(G)\leq N$. 
\item If $L'\leq (H\times G)$, then $[(H\times G)/L']\times_G\varphi_N^G=0$ unless $p_2(L')=G$ and $k_2(L')\cap\Phi(G)\leq\nolinebreak N$.
\end{enumerate}
\end{enumerate}
\end{mth}
\pf For $N\normal G$, set $u_N^G=Y_{\Delta_N(G)}$. Since $\Delta_N(G)*\Delta_M(G)=\Delta_{NM}(G)$ for any normal subgroups $N$ and $M$ of $G$, it follows from Corollary~\ref{sub Frattini} that if either $N$ or $M$ is contained in $\Phi(G)$, then $u_N^Gu_M^G=u_{NM}^G$.\par
Now Assertion 1 follows from the following general procedure for building orthogonal idempotents (see \cite{essential} Theorem 10.1 for details): we have a finite lattice $P$ (here $P$ is the lattice of normal subgroups of $G$ contained in $\Phi(G)$), and a set of elements $g_x$ of a ring $A$, for $x\in P$ (here $A=\mathcal{E}(G)$ and $g_N=u_N^G$), with the property that $g_xg_y=g_{x\vee y}$ for any $x,y\in P$, and $g_0=1$, where 0 is the smallest element of $P$ (here this element is the trivial subgroup of $G$, and $u_\un^G=Y_{\Delta_\un(G)}=\widetilde{e_G^G}$). The the elements $f_x$ defined for $x\in P$ by
$$f_x=\sumb{y\in P}{x\leq y}\mu(x,y)g_y\mvirg$$
where $\mu$ is the M\"obius function of $P$, are orthogonal idempotents of $A$, and their sum is equal to the identity element of $A$. This is exactly Assertion 1 (since $f_x=\varphi_N^G$ here, for $x=N\in P$).\par
Let $L\leq (G\times H)$, then by~\ref{factorize}
$$\varphi_N^G\times_G[(G\times H)/L]=\varphi_N^G\circ\Ind_{p_1(L)}^G\circ[(p_1(L)\times H)/L]=0$$
unless $p_1(L)=G$, by Corollary~\ref{essentiel frattini}. And if $p_1(L)=G$, then by~\ref{factorize}
$$\varphi_N^G\times_G[(G\times H)/L]=\varphi_N^G\circ\Inf_{G/k_1(L)}^G\circ[(G/k_1(L)\times H)/L_1\mvirg$$
for some subgroup $L_1$ of $(G/k_1(L)\times H)$. Again, by Corollary~\ref{essentiel frattini} this is equal to 0 unless $k_1(L)\cap\Phi(G)\leq N$.
The proof of Assertion~(b) is similar. Alternatively, one can take opposite bisets in (a).\endpf
\begin{mth}{Proposition} \label{phiYL}Let $G$ be a finite group with order invertible in~$\Cf$.
\begin{enumerate}
\item Let $L\in\Sigma(G,G)$. Then 
$$\varphi_\un^GY_L=\sumb{N\normal G}{N\leq\Phi(G)}\mu_{\normal G}(\un,N)Y_{(N\times \un)L}\mpoint$$
This is non zero if and only if $k_1(L)\cap \Phi(G)=\un$. 
Similarly 
$$Y_L\varphi_\un^G=\sumb{N\normal G}{N\leq\Phi(G)}\mu_{\normal G}(\un,N)Y_{L(\un\times N)}\mvirg$$
and $Y_L\varphi_1^G\neq 0$ if and only if $k_2(L)\cap\Phi(G)=\un$.
\item The elements $\varphi_\un^GY_L$ (resp. $Y_L\varphi_\un^G$), when $L$ runs through a set of representatives of conjugacy classes of elements of $\Sigma(G,G)$ such that $k_1(L)\cap \Phi(G)=\un$ (resp $k_2(L)\cap\Phi(G)=\un$), form an $\Cf$-basis of the right ideal $\varphi_\un^G\mathcal{E}(G)$ (resp. the left ideal $\mathcal{E}(G)\varphi_\un^G$) of $\mathcal{E}(G)$.
\end{enumerate}
\end{mth}
\pf Let $L\in\Sigma(G,G)$. By Proposition~\ref{idempotents E}, we have
\begin{eqnarray*}
\varphi_\un^GY_L\!\!\!&=&\!\!\!\widetilde{e_G^G}\times_G\Big(\!\sumb{N\normal G}{N\leq \Phi(G)}\!\mu_{\normal G}(\un,N)[(G\times G)/\Delta_{N}(G)]\Big)\times_G[(G\times G)/L]\times_G\widetilde{e_G^G}\\
&=&\widetilde{e_G^G}\times_G\Big(\sumb{N\normal G}{N\leq \Phi(G)}\mu_{\normal G}(\un,N)[(G\times G)/\Delta_N(G)*L]\Big)\times_G\widetilde{e_G^G}\\
&=&\widetilde{e_G^G}\times_G\Big(\sumb{N\normal G}{N\leq \Phi(G)}\mu_{\normal G}(\un,N)[(G\times G)/(N\times\un)L]\Big)\times_G\widetilde{e_G^G}\mpoint\\
&=&\sumb{N\normal G}{N\leq \Phi(G)}\mu_{\normal G}(\un,N)Y_{(N\times\un)L}\mpoint\\
\end{eqnarray*}
Set $M=k_1(L)\cap\Phi(G)$. Then $M\normal G$, and $(N\times\un)L=(NM\times\un)L$ for any normal subgroup $N$ of $G$ contained in $\varphi(G)$. Thus
\begin{equation}\label{produit}
\varphi_\un^GY_L=\sumb{P\normal G}{M\leq P\leq\Phi(G)}\Big(\sumb{N\normal G}{NM=P}\mu_{\normal G}(\un,N)\Big)Y_{(P\times\un)L}\mpoint
\end{equation}
If $M\neq \un$, then $\Big(\dsp{\sumb{N\normal G}{NM=P}\mu_{\normal G}(\un,N)\Big)=0}$ for any $P\normal G$ with $M\leq P\leq \Phi(G)$. Hence $\varphi_\un^GY_L=0$ in this case. And if $M=\un$, Equation (\ref{produit}) reads
$$\varphi_\un^GY_L=\sumb{P\normal G}{P\leq\Phi(G)}\mu_{\normal G}(\un,P)Y_{(P\times \un)L}\mpoint$$
The element $Y_{(P\times \un)L}$ is equal to $Y_L$ if and only if $(P\times \un)L$ is conjugate to $L$. This implies that $k_1\big((P\times \un)L\big)$ is conjugate to (hence equal to) $k_1(L)$. Thus $P\leq k_1\big((P\times \un)L\big)\leq k_1(L)\cap\Phi(G)$, hence $P=\un$. So the coefficient of $Y_L$ in $\varphi_\un^GY_L$ is equal to 1, hence $\varphi_\un^GY_L\neq 0$. The remaining statements of Assertion 1 follow by taking opposite bisets.\par
Assertion 2 follows from Proposition~\ref{product of YL}, and from the fact that the coefficient of $Y_L$ in $\varphi_\un^GY_L$ is equal to 1 when $k_1(L)\cap \Phi(G)=\un$.\endpf
\begin{mth}{Corollary} Let $G$ be a finite group of order invertible in $\Cf$. If every minimal (non-trivial) normal subgroup of $G$ is contained in $\Phi(G)$, then $\varphi_\un^G$ is central in $\mathcal{E}(G)$, and the algebra $\varphi_\un^G\CE(G)$ is isomorphic to $\Cf \Out(G)$.
\end{mth}
\pf Indeed if $L\in\Sigma(L,L)$ and $\varphi_\un^GY_L\neq 0$, then $k_1(L)\cap\Phi(G)=\un$. It follows that $k_1(L)$ contains no minimal normal subgroup of $G$, and then $k_1(L)=\un$. Equivalently $q(L)\cong p_1(L)/k_1(L)\cong G\cong p_2(L)/k_2(L)$, i.e. $k_2(L)=G$ also, or equivalently $k_2(L)\cap\Phi(G)=\un$. Hence $\varphi_\un^GY_L\neq 0$ if and only if $Y_L\varphi_\un^G\neq 0$, and in this case, there exists an automorphism $\theta$ of $G$ such that 
$$L=\Delta_\theta(G)=\{\big(\theta(x),x\big)\mid x\in G\}\mpoint$$
In this case for any normal subgroup $N$ of $G$ contained in $\Phi(G)$
\begin{eqnarray*}
(N\times \un)L&=&\{(a,b)\in G\times G\mid a\theta(b)^{-1}\in N\}\\
&=&\{(a,b)\in G\times G\mid a^{-1}\theta(b)\in N\}\\
&=&L\big(\un\times \theta^{-1}(N)\big)\mpoint
\end{eqnarray*}
Now $N\mapsto\theta^{-1}(N)$ is a permutation of the set of normal subgroups of $G$ contained in $\Phi(G)$. Moreover $\mu_{\normal G}(\un,N)=\mu_{\normal G}\big(\un,\theta^{-1}(N)\big)$. \par
It follows that $\varphi_\un^GY_L=Y_L\varphi_\un^G$, so $\varphi_\un^G$ is central in $\CE(G)$. Moreover the map $\theta\in\Aut(G)\mapsto \varphi_\un^GY_{\Delta_\theta(G)}$ clearly induces an algebra isomorphism $\Cf \Out(G)\to \varphi_\un^G\CE(G)$.\findemo 
\begin{mth}{Theorem} \label{central}Let $G$ be a finite group with order invertible in~$\Cf$. If $G$ is nilpotent, then $\varphi_\un^G$ is a central idempotent of $\mathcal{E}(G)$.
\end{mth}
\pf Let $L\in\Sigma(G,G)$. Setting $Q=q(L)$, there are two surjective group homomorphisms $s,t:G\to Q$ such that $L=\{(x,y)\in G\times G\mid s(x)=t(y)\}$. Then $k_1(L)=\Ker\,s$ and $k_2(L)=\Ker\,t$. Now by Proposition~\ref{phiYL}
$$\varphi_\un^GY_L=\sumb{N\normal G}{N\leq\Phi(G)}\mu_{\normal G}(\un,N)Y_{(N\times \un)L}\mvirg$$
and this is non zero if and only if $\Ker\,s\cap\Phi(G)=\un$. 
Now $s\big(\Phi(G)\big)$ is equal to $\Phi(Q)$ since $G$ is nilpotent: indeed $G=\prod_pG_p$ (resp. $Q=\prod_pQ_p$) is the direct product of its $p$-Sylow subgroups $G_p$ (resp.~$Q_p$), and $s$ induces a surjective group homomorphism $G_p\to Q_p$, for any prime~$p$. Moreover $\Phi(G)=\prod_p\Phi(G_p)$ (resp. $\Phi(Q)=\prod_p\Phi(Q_p)$). Finally $\Phi(G_p)$ is the subgroup of $G_p$ generated by commutators and $p$-powers of elements of $G_p$, hence it maps by $s$ {\em onto} the subgroup of $Q_p$ generated by commutators and $p$-powers of elements of $Q_p$, that is $\Phi(Q_p)$. Similarly $t\big(\Phi(G)\big)=\Phi(Q)$. \par
If $\Ker\,s\cap\Phi(G)=\un$, it follows that $s$ induces an isomorphism from $\Phi(G)$ to $\Phi(Q)$. Then the surjective homomorphism $\Phi(G)\to \Phi(Q)$ induced by $t$ is also an isomorphism, and in particular $\Ker\,t\cap\Phi(G)=\un$.\par
Let $D=L\cap\big(\Phi(G)\times\Phi(G)\big)$. Then $k_1(D)\subseteq k_1(L)\cap \Phi(G)=\Ker\,s\cap\Phi(G)$, hence $k_1(D)=\un$. Similarly $k_2(L)\subseteq k_2(L)\cap \Phi(G)=\Ker\,t\cap\Phi(G)=\un$, hence $k_2(D)=\un$. Since $s\big(\Phi(G)\big)=\Phi(Q)=t\big(\Phi(G)\big)$, we have moreover $p_1(D)=\Phi(G)=p_2(D)$. It follows that there is an automorphism $\alpha$ of $\Phi(G)$ such that $D=\{\big(x,\alpha(x)\big)\mid x\in\Phi(G)\}$. \par
Moreover for any $y\in G$, there exists $z\in G$ such that $(y,z)\in L$. It follows that $\big(x^y,\alpha(x)^z\big)\in D$ for any $x\in\Phi(G)$, that is $\alpha(x^y)=\alpha(x)^z$. In particular if $N$ is a normal subgroup of $G$ contained in $\Phi(G)$, then so is $\alpha(N)$. Hence $\alpha$ induces an automorphism of the poset of normal subgroups of $G$ contained in $\Phi(G)$. In particular $\mu_{\normal G}(\un,N)=\mu_{\normal G}\big(\un,\alpha(N)\big)$. \par
Moreover for $n\in N$ and $(y,z)\in L$, we have
$$(n,1)(y,z)=(y,z)(n^y,1)=(y,z)\big(n^y,\alpha(n^y)\big)\big(1,\alpha(n^y)^{-1}\big)\mpoint$$
Since $\big(n^y,\alpha(n^y)\big)\in D\leq L$, we have $(N\times \un)L=L\big(\un\times \alpha(N)\big)$. It follows that
\begin{eqnarray*}
\varphi_\un^GY_L&=& \sumb{N\normal G}{N\leq \Phi(G)}\mu_{\normal G}(\un,N)Y_{(N\times \un)L}=\sumb{N\normal G}{N\leq \Phi(G)}\mu_{\normal G}(\un,N)Y_{L(\un\times\alpha(N))}\\
&=&\sumb{N\normal G}{N\leq \Phi(G)}\mu_{\normal G}\big(\un,\alpha(N)\big)Y_{L(\un\times\alpha(N))}=\sumb{N\normal G}{N\leq \Phi(G)}\mu_{\normal G}\big(\un,N)Y_{L(\un\times N)}\\
&=&Y_L\varphi_\un^G\mvirg
\end{eqnarray*}
as was to be shown.\endpf
\begin{rem}{Remark} When $G$ is not nilpotent, it is not true in general that $\varphi_\un^G$ is central in $\mathcal{E}(G)$. This is because $t\big(\Phi(G)\big)$ need not be equal to $\Phi(Q)$ for a surjective group homomorphism $t:G\to Q$. For example, there is a surjection $t$ from the group $G=C_4\times(C_5\rtimes C_4)$ to $Q=C_4$ with kernel $C_4\times C_5$ {\em containing} $\Phi(G)=C_2\times\un$, and another surjection $s:G\to Q$ with kernel $\un\times(C_5\rtimes C_4)$ intersecting trivially $\Phi(G)$. In this case, the group $L=\{(x,y)\in G\times G\mid s(x)=t(y)\}$ is in $\Sigma(G,G)$, and $k_1(L)\cap \Phi(G)=\un$, but $k_2(L)\cap\Phi(G)=\Phi(G)\neq\un$. By Proposition~\ref{phiYL}, we have $\varphi_\un^GY_L\neq 0$ and $Y_L\varphi_\un^G=0$, so $\varphi_\un^G$ is not central in $\mathcal{E}(G)$.
\end{rem} 
\section{Idempotents in $\Cf B(G,G)$}
\begin{mth}{Definition} When $G$ is a finite group, a {\em section} $(T,S)$ of $G$ is a pair of subgroups of $G$ such that $S\normal T$.\par
A section $(T,S)$ is called {\em minimal} (cf. \cite{vanishing}) if $S\leq \Phi(T)$. Let $\mathcal{M}(G)$ denote the set of minimal sections of $G$.\par
A group $H$ is called a subquotient of $G$ (notation $H\sqsubseteq G$) if there exists a section $(T,S)$ of $G$ such that $T/S\cong H$.
\end{mth}
A section $(T,S)$ is minimal if and only if the only subgroup $H$ of~$T$ such that $H/(H\cap S)\cong T/S$ is $T$ itself.
\begin{mth}{Notation} Let $G$ be a finite group, and let $(T,S)$ be a section of $G$. 
\begin{enumerate}
\item Let $\Indinf_{T/S}^G\in B(G,T/S)$ denote (the isomorphism class of) the $(G,T/S)$-biset $G/S$, and let $\Defres_{T/S}^G\in B(T/S,G)$ denote (the isomorphism class of) the $(T/S,G)$-biset $S\dom G$.
\item Let $\Cf$ be a commutative ring in which the order of $G$ is invertible. Let $u_{T,S}^G\in \Cf B(G,T/S)$ be defined by
$$u_{T,S}^G=\Indinf_{T/S}^G\varphi_\un^{T/S}\mvirg$$
and let $v_{T,S}^G\in \Cf B(T/S,G)$ be defined by
$$v_{T,S}^G=\varphi_\un^{T/S}\Defres_{T/S}^G\mpoint$$
\end{enumerate}
\end{mth}
\begin{rem}{Remark} Observe that $v_{T,S}^G=(u_{T,S}^G)\op$: indeed $(G/S)\op \cong S\dom G$, and $(\varphi_\un^{T/S})\op=\varphi_\un^{T/S}$.
\end{rem}
\begin{mth}{Theorem} \label{u et v}Let $G$ be a finite group with order invertible in~$\Cf$.\begin{enumerate}
\item If $(T,S)$ and $(T',S')$ are minimal sections of $G$, then
$$v_{T',S'}^Gu_{T,S}^G=0$$
unless $(T,S)$ and $(T',S')$ are conjugate in $G$.
\item If $(T,S)$ is a minimal section of $G$, then
$$v_{T,S}^Gu_{T,S}^G=\varphi_\un^{T/S}\Big(\sum_{g\in N_G(T,S)/T} \Iso(c_g)\Big)\mvirg$$
where $N_G(T,S)=N_G(T)\cap N_G(S)$, and $c_g$ is the automorphism of $T/S$ induced by conjugation by $g$.
\end{enumerate}
\end{mth}
\pf Indeed $(S'\dom G)\times_G(G/S)\cong S'\dom G/S$ as a $(T'/S',T/S)$-biset. Hence
$$v_{T',S'}^Gu_{T,S}^G=\varphi_\un^{T'/S'}\Big(\sum_{g\in T'\dom G/T} S'\dom T'gT/S\Big)\varphi_\un^{T/S}\mpoint$$
For any $g\in G$, the $(T'/S',T/S)$-biset $S'\dom T'gT/S$ is transitive, isomorphic to $\big((T'/S')\times (T/S)\big)/L_g$, where
$$L_g=\{(t'S',tS)\in (T'/S')\times(T/S)\mid t'gt^{-1}\in S'gS\}\mpoint$$
Then $t'S'\in p_1(L_g)$ if and only if $t'\in S'\cdot gTg^{-1}\cap T'$. Hence 
$$p_1(L_g)=({^gT}\cap T')S'/S'\mpoint$$
Similarly $p_2(L_g)=(T'^g\cap T)S/S$. In particular $p_1(L_g)=T'/S'$ if and only if $({^gT}\cap T')S'=T'$, i.e. ${^gT}\cap T'=T'$, since $S'\leq\Phi(T')$. Thus 
$p_1(L_g)=T'/S'$ if and only if $T'\leq {^gT}$. Similarly $p_2(L_g)=T/S$ if and only if $T\leq T'^g$. By Theorem~\ref{idempotents E}, it follows that $\varphi_\un^{T'/S'}(S'\dom T'gT/S)\varphi_\un^{T/S}=0$ unless $T'={^gT}$.\par
Assume now that $T'={^gT}$. Then $t'S'\in k_1(L_G)$ if and only if $t'$ lies in  $S'\cdot gSg^{-1}\cap T'$. Hence 
$$k_1(L_g)=({^gS}\cap T')S'/S'\mvirg$$
and similarly $k_2(L_g)=(S'^g\cap T)S/S$. But since $S\leq \Phi(T)$ and $S\normal T$, it follows that $^gS\normal {^gT}=T'$ and $^gS\leq{^g\Phi(T)}=\Phi(T')$. Hence $^gS\cdot S'/S'$ is contained in $k_1(L_g)\cap\Phi(T')/S'$. Moreover $\Phi(T')/S'=\Phi(T'/S')$, as 
$$\Phi(T'/S')=\mathop{\bigcap}_{S'\leq M'< T'}\limits(M'/S')=\mathop{\bigcap}_{M'< T'}\limits(M'/S')=(\mathop{\bigcap}_{M'<T'}\limits M')/S'=\Phi(T')/S'\mvirg$$
where $M'$ runs through maximal subgroups of $T'$, which all contain $S'$ since $S'\leq \Phi(T')$.\par
It follows that if $k_1(L_g)\cap \Phi(T'/S')=\un$, then $^gS\cdot S'=S'$, that is $^gS\leq S'$. Similarly if $k_2(L_g)\cap \Phi(T/S)=\un$, then $S'^g\leq S$. By Theorem~\ref{idempotents E}, it follows that $\varphi_\un^{T'/S'}(S'\dom T'gT/S)\varphi_\un^{T/S}=0$ unless $T'={^gT}$ and $S'={^gS}$. This proves Assertion 1.\mpn
For Assertion 2, the same computation shows that 
$$v_{T,S}^Gu_{T,S}^G=\sum_{g\in N_G(T,S)/T}\varphi_\un^{T/S}(S\dom TgT/S)\varphi_\un^{T/S}\mpoint$$
But $S\dom TgT/S=gT/S$ if $g\in N_G(T,S)$, and this $(T/S,T/S)$-biset is isomorphic to $\Iso(c_g)$. Assertion 2 follows, since moreover $\varphi_\un^{T/S}$ commutes with any biset of the form $\Iso(\theta)$, where $\theta$ is an automorphism of $T/S$.\endpf
\begin{mth}{Notation} For a minimal section $(T,S)$ of the group $G$, set
$$\epsilon_{T,S}^G=\fracb{1}{|N_G(T,S):T|}u_{T,S}^Gv_{T,S}^G=\fracb{1}{|N_G(T,S):T|}\Indinf_{T/S}^G\varphi_\un^G\Defres_{T/S}^G\in\Cf B(G,G)\mpoint$$
\end{mth}
Note that $\epsilon_{T,S}^G=\epsilon_{^gT,^gS}^G$ for any $g\in G$, and that $\epsilon_{G,N}^G=\varphi_N^G$ when $N\normal G$ and $N\leq\Phi(G)$, by Proposition~\ref{phi1}.
\begin{mth}{Proposition} \label{formule epsilon}Let $(T,S)$ be a minimal section of $G$. Then
$$\epsilon_{T,S}^G\!=\frac{1}{|N_G(T,S)|}\!\!\!\!\!\!\sumb{X\leq T, M\normal T}{\rule{0ex}{1.6ex}S\leq M\leq\Phi(T)\leq X\leq T}|X|\mu(X,T)\mu_{\normal T}(S,M)\,\Indinf_{X/M}^G\circ\Defres_{X/M}^G\mpoint$$
\end{mth}
\pf This is a straightforward consequence of the above definition of $\epsilon_{T,S}^G$, and from Assertion~2 of Proposition~\ref{phinG}.\endpf
\begin{mth}{Theorem} Let $G$ be a finite group with order invertible in~$\Cf$, let $[\mathcal{M}(G)]$ be a set of representatives of conjugacy classes of minimal sections of~$G$. Then the elements $\epsilon_{T,S}^G$, for $(T,S)\in[\mathcal{M}(G)]$, are orthogonal idempotents of $\Cf B(G,G)$, and their sum is equal to the identity element of $\Cf B(G,G)$.
\end{mth}
\pf Let $(T,S)$ and $(T',S')$ be distinct elements of $[\mathcal{M}(G)]$. Then
$$\epsilon_{T',S'}^G\epsilon_{T,S}^G=\fracb{1}{|N_G(T',S'):T'|}\fracb{1}{|N_G(T,S):T|}u_{T',S'}^Gv_{T',S'}^Gu_{T,S}^Gv_{T,S}^G=0\mvirg$$
since $v_{T',S'}^Gu_{T,S}^G=0$ by Theorem~\ref{u et v}. Moreover:
\begin{eqnarray*}
\sum_{(T,S)\in [\mathcal{M}(G)]}\epsilon_{T,S}^G&=&\sum_{(T,S)\in [\mathcal{M}(G)]}\fracb{1}{|N_G(T,S):T|}u_{T,S}^Gv_{T,S}^G\\
&=&\sum_{(T,S)\in \mathcal{M}(G)}\fracb{1}{|G:T|}u_{T,S}^Gv_{T,S}^G\\
&=&\sum_{(T,S)\in \mathcal{M}(G)}\fracb{1}{|G:T|}\Indinf_{T/S}^G\varphi_\un^{T/S}\Defres_{T/S}^G
\end{eqnarray*}
Now $\varphi_\un^{T/S}=\widetilde{e_{T/S}^{T/S}}f_{T/S}$ by Proposition~\ref{phinG}, where
$$f_{T/S}=\sumb{N/S\normal (T/S)}{N/S\leq \Phi(T/S)}\mu_{\normal G}(\un,N/S)[\big((T/S)\times(T/S)\big)/\Delta_{N/S}(T/S)]\mpoint$$
Hence $\varphi_\un^{T/S}=\widetilde{e_{T/S}^{T/S}}\Def_{T/S}^T\Inf_{T/S}^Tf_{T/S}$, and
$$\sum_{(T,S)\in [\mathcal{M}(G)]}\epsilon_{T,S}^G=\sum_{(T,S)\in \mathcal{M}(G)}\fracb{1}{|G:T|}\Ind_T^G\Inf_{T/S}^T\widetilde{e_{T/S}^{T/S}}\Def_{T/S}^T\Inf_{T/S}^Tf_{T/S}\Def_{T/S}^T\Res_T^G\mpoint$$
Now $\Inf_{T/S}^T\widetilde{e_{T/S}^{T/S}}\Def_{T/S}^T=\widetilde{\Inf_{T/S}^Te_{T/S}^{T/S}}$, and $\Inf_{T/S}^Te_{T/S}^{T/S}$ is equal to the sum over subgroups $X$ or $T$ such that $XS=T$, up to conjugation, of the idempotents~$e_X^T$. Since $S\leq\Phi(T)$, the only subgroup $X$ of $T$ such that $XS=T$ is $T$ itself. Hence
$$\Inf_{T/S}^T\widetilde{e_{T/S}^{T/S}}\Def_{T/S}^T=\widetilde{e_T^T}\mpoint$$
On the other hand
$$\Inf_{T/S}^T[\big((T/S)\times(T/S)\big)/\Delta_{N/S}(T/S)]\Def_{T/S}^T=[(T\times T)/\Delta_N(T)]\mpoint$$
It follows that the sum $\Sigma=\sum_{(T,S)\in [\mathcal{M}(G)]}\limits\epsilon_{T,S}^G$ is equal to
\begin{eqnarray*}
\Sigma&=&\sum_{(T,S)\in \mathcal{M}(G)}\fracb{1}{|G:T|}\Ind_T^G\widetilde{e_T^T}\sumb{N\normal T}{S\leq N\leq \Phi(T)}\mu_{\normal T}(S,N)[(T\times T)/\Delta_N(T)]\Res_T^G\\
&=&\sum_{(T,S)\in \mathcal{M}(G)}\fracb{1}{|G:T|}\Ind_T^G\widetilde{e_T^T}\varphi_{S}^T\Res_T^G\;\;\hbox{[by definition of $\varphi_S^T$]}\\
&=&\sum_{T\leq G}\fracb{1}{|G:T|}\Ind_T^G\widetilde{e_T^T}\sumb{S\normal T}{S\leq\Phi(T)}\varphi_{S}^T\Res_T^G\\
&=&\sum_{T\leq G}\fracb{1}{|G:T|}\Ind_T^G\widetilde{e_T^T}\Res_T^G\;\;\hbox{[by Theorem~\ref{idempotents E}]}\\
&=&\sum_{T\leq G}\fracb{1}{|G:T|}\widetilde{\Ind_T^Ge_T^T}\;\;\hbox{[by Lemma~\ref{tilde compatibility}]}\\
&=&\sum_{T\leq G}\fracb{1}{|G:N_G(T)|}\widetilde{e_T^G}\;\;\hbox{[by (\ref{ehg})]}\\
&=&\sum_{T\in[s_G]}\widetilde{e_T^G}=\widetilde{G/G}=[(G\times G)/\Delta(G)]\mpoint
\end{eqnarray*}
So the sum $\Sigma$ is equal to the identity of $\Cf B(G,G)$. Since $\epsilon_{T,S}^G\epsilon_{T',S'}^G=0$ if $(T,S)$ and $(T',S')$ are distinct elements of $[\mathcal{M}(G)]$, it follows that for any $(T,S)\in [\mathcal{M}(G)]$
$$\epsilon_{T,S}^G=\epsilon_{T,S}^G\Sigma=(\epsilon_{T,S}^G)^2\mvirg$$
which completes the proof of the theorem.\endpf
\section{Application to biset functors}
\begin{mth}{Notation} Let $F$ be a biset functor over $\Cf$. When $G$ is a finite group with order invertible in $\Cf$, we set
$$\delta_\Phi F(G)=\varphi_\un^GF(G)$$
\end{mth}
\begin{mth}{Proposition} Let $F$ be a biset functor over $\Cf$. Then for any finite group $G$ with order invertible in~$\Cf$, the $\Cf$-submodule $\delta_\Phi F(G)$ of $F(G)$ is the set of elements $u\in F(G)$ such that
$$\left\{\begin{array}{ll} \Res_H^Gu=0&\forall H<G\\\Def_{G/N}^Gu=0&\forall N\normal G,\;N\cap\Phi(G)\neq\un\end{array}\right.\mpoint$$
\end{mth}
\pf If $u\in \delta_\Phi F(G)=\varphi_\un^GF(G)$, then $\Res_H^Gu=0$ for any proper subgroup $H$ of $G$, and $\Def_{G/N}^Gu=0$ for any $N\normal G$ such that $N\cap\Phi(G)\neq\un$, by Corollary~\ref{essentiel frattini}. \par
Conversely, if $u\in F(G)$ fulfills the two conditions of the proposition, then $\widetilde{e_G^G}u=u$, because $\widetilde{e_G^G}$ is equal to the identity element $[(G\times G)/\Delta(G)]$ of $\Cf B(G,G)$, plus a linear combination of elements of the form $[(G\times G)/\Delta(H)]=\Ind_H^G\circ\Res_H^G$, for proper subgroups $H$ of $G$. Similarly $\Inf_{G/N}^G\Def_{G/N}^Gu=0$ for any non-trivial normal subgroup of $G$ contained in $\Phi(G)$, thus $\varphi_\un^Gu=u$.\endpf
\begin{rem}{Remark} Since $\Def_{G/N}^G=\Def_{G/N}^{G/M}\circ\Def_{G/M}^G$, where $M=N\cap\Phi(G)$, saying that $\Def_{G/N}^Gu=0$ for any $N\normal G$ with $N\cap\Phi(G)\neq\un$ is equivalent to saying that $\Def_{G/N}^Gu=0$ for any non trivial normal subgroup $N$ of $G$ contained in $\Phi(G)$.
\end{rem}
\pagebreak[3]
\begin{mth}{Theorem}\label{evaluation} Let $F$ be a biset functor over $\Cf$. Then for any finite group~$G$ with order invertible in~$\Cf$, the maps
$$\xymatrix@R=1ex{
F(G)\ar[r]<.5ex>&*!U(.3){\dirsum{(T,S)\in [\mathcal{M}(G)]} \big(\delta_\Phi F(T/S)\big)^{N_G(T,S)/T}}\ar[l]<.5ex>\\
\hspace{4ex}w\ar@{|->}[r]^-V&*!U(.4){\dirsum{(T,S)} \frac{1}{|N_G(T,S):T|}v_{T,S}^Gw\hspace{5.4ex}}\\
*!U(.4){\sum_{(T,S)}\limits u_{T,S}^Gw_{T,S}\;}&*!U(.6){\dirsum{(T,S)}w_{T,S}\hspace{15.5ex}}\ar@{|->}[l]_-U
}
$$
are well defined isomorphisms of $\Cf$-modules, inverse to one other.
\end{mth}
\pf We have first to check that if $w\in F(G)$, then the element $v_{T,S}^Gw$ of $\varphi_\un^{T/S}F(T/S)=\delta_\Phi F(T/S)$ is invariant under the action of $N_G(T,S)/T$. But for any $g\in N_G(T/S)$
$$\Iso(c_g)v_{T,S}^G=v_{^gT,^gS}^G\Iso(c_g)=v_{T,S}^G\Iso(c_g)\mvirg$$
where $\Iso(c_g):F(G)\to F(G)$ on the right hand side is conjugation by $g$, that is an inner automorphism, hence the identity map, for $g\in G$.\par
Now for $w\in F(G)$
\begin{eqnarray*}
UV(w)&=&\sum_{(T,S)\in [\CM(G)]}\fracb{1}{|N_G(T,S):T|}u_{T,S}^Gv_{T,S}^Gw\\
&=&\sum_{(T,S)\in[\CM(G)]}\epsilon_{T,S}^Gw=w\mvirg
\end{eqnarray*}
so $UV$ is the identity map of $F(G)$.\par
Conversely, if $w_{T,S}\in \big(\delta_\Phi F(T/S)\big)^{N_G(T,S)/T}$, for $(T,S)\in[\CM(G)]$, then
\begin{eqnarray*}
VU\big(\dirsum{(T,S)\in[\CM(G)]}w_{T,S}\big)\!\!\!&=&\!\!\dirsum{(T,S)\in[\CM(G)]}\!\sum_{(T',S')\in[\CM(G)]}\fracb{1}{|N_G(T,S):T|}v_{T,S}^Gu_{T',S'}^Gw_{T',S'}\\
\!\!\!&=&\!\!\!\dirsum{(T,S)\in[\CM(G)]}\!\fracb{1}{|N_G(T,S):T|}v_{T,S}^Gu_{T,S}^Gw_{T,S}\\
\!\!\!&=&\!\!\!\dirsum{(T,S)\in[\CM(G)]}\!\fracb{1}{|N_G(T,S):T|}\sum_{g\in N_G(T,S)/T}\Iso(c_g)w_{T,S}\\
\!\!\!&=&\!\!\!\dirsum{(T,S)\in[\CM(G)]}w_{T,S}\mvirg\\
\end{eqnarray*}
so $VU$ is also equal to the identity map.\endpf
\section{Atoric $p$-groups}
For the remainder of the paper, we denote by $p$ a (fixed) prime number.
\begin{mth}{Notation and Definition} \begin{itemize}
\item If $P$ is a finite $p$-group, let $\Omega_1P$ denote the subgroup of $P$ generated by the elements of order $p$.
\item A finite $p$-group $P$ is called {\em atoric} if it does not admit any decomposition $P=E\times Q$, where $E$ is a non-trivial elementary abelian $p$-group. Let $\At_p$ denote the class of atoric $p$-groups, and let $[\At_p]$ denote a set of representatives of isomorphism classes in $\At_p$.
\end{itemize}
\end{mth}
The terminology ``atoric" is inspired by~\cite{quillen}, where elementary abelian $p$-groups are called {\em $p$-tori}. Atoric $p$-groups have been considered (without naming them) in~\cite{ideal}, Example 5.8.
\begin{mth}{Lemma} \label{NcapPhi}Let $P$ be a finite $p$-group, and $N$ be a normal subgroup of~$P$.  The following conditions are equivalent:
\begin{enumerate}
\item $N\cap \Phi(P)=\un$
\item $N$ is elementary abelian and central in $P$, and admits a complement in~$P$.
\item $N$ is elementary abelian and there exists a subgroup $Q$ of $P$ such that $P=N\times Q$.
\end{enumerate}
\end{mth}
\pf\mpn
\fbox{$1\Rightarrow 3$} Let $N\normal P$ with $N\cap\Phi(P)=\un$. Then $N$ maps injectively in the elementary abelian $p$-group $P/\Phi(P)$, so $N$ is elementary abelian.  Let $Q/\Phi(P)$ be a complement of $N\Phi(P)/\Phi(P)$ in $P/\Phi(P)$. Then $Q\geq \Phi(P)\geq [P,P]$, so $Q$ is normal in $P$. Moreover $Q\cdot N=P$ and $Q\cap N\Phi(P)=(Q\cap N)\Phi(P)=\Phi(P)$, thus $Q\cap N\leq \Phi(P)\cap N=\un$. Now $N$ and $Q$ are normal subgroups of $P$ which intersect trivially, hence they centralize each other. It follows that $P=N\times Q$.\mpn
\fbox{$3\Rightarrow 2$} This is clear.\mpn
\fbox{$2\Rightarrow 1$} If $P=N\cdot Q$ for some subgroup $Q$ of $P$, and if $N$ is central in $P$, then $P=N\times Q$. Thus $\Phi(P)=\un\times \Phi(Q)$, as $N$ is elementary abelian. Then $N\cap \Phi(P)\leq N\cap Q=\un$.\endpf
\begin{mth}{Lemma} \label{atoric}Let $P$ be a finite $p$-group. The following conditions are equivalent:
\begin{enumerate}
\item $P$ is atoric.
\item If $N\normal P$ and $N\cap\Phi(P)=\un$, then $N=\un$.
\item $\Omega_1Z(P)\leq \Phi(P)$.
\end{enumerate}
\end{mth}
\pf\mpn
\fbox{$1\Rightarrow 2$} Suppose that $P$ is atoric. Let $N\normal P$ with $N\cap\Phi(P)=\un$. Then by Lemma~\ref{NcapPhi}, the group $N$ is elementary abelian and there exists a subgroup $Q$ of $P$ such that $P=N\times Q$. Hence $N=\un$.\mpn
\fbox{$2\Rightarrow 3$} Suppose now that Assertion~2 holds. If $x$ is a central element of order~$p$ of~$P$, then the subgroup $N$ of $P$ generated by $x$ is normal in $P$, and non trivial. Then $N\cap \Phi(P)\neq \un$, hence $N\leq \Phi(P)$ since $N$ has order $p$, thus $x\in\Phi(P)$. \mpn
\fbox{$3\Rightarrow 1$} Finally, if Assertion~3 holds, and if $P=E\times Q$ for some subgroups $E$ and $Q$ of $P$ with $E$ elementary abelian, then $\Phi(P)=\un\times\Phi(Q)$. Moreover $E\leq\Omega_1Z(P)\leq\Phi(P)\leq Q$, so $E=E\cap Q=\un$, and $P$ is atoric.\endpf
\pagebreak[3]
\begin{mth}{Proposition} \label{Pzero}Let $P$ be a finite $p$-group, and $N$ be a maximal normal subgroup of $P$ such that $N\cap\Phi(P)=\un$. Then: 
\begin{enumerate}
\item The group $N$ is elementary abelian and there exists a subgroup $T$ of $P$ such that $P=N\times T$.
\item The group $P/N\cong T$ is atoric.
\item If $Q$ is an atoric $p$-group and $s:P\twoheadrightarrow Q$ is a surjective group homomorphism, then $s(T)=Q$. In particular $Q$ is isomorphic to a quotient of $T$.
\end{enumerate}
\end{mth}
\pf (1) This follows from Lemma~\ref{NcapPhi}.\mpn
(2) By (1), there exists $T\leq P$ such that $P=N\times T$. In particular $P/N\cong T$. Now if $T=E\times S$, for some subgroups $E$ and $S$ of $T$ with $E$ elementary abelian, then $P\cong P_1=(N\times E)\times S$, and $N\times E$ is an elementary abelian normal subgroup of $P_1$ which intersects trivially $\Phi(P_1)=\Phi(S)$. By maximality of $N$, it follows that $E=\un$, so $T\cong P/N$ is atoric.\mpn
(3) Let $s:P\twoheadrightarrow Q$ be a surjective group homomorphism, where $Q$ is atoric. By (1), the group $N$ is elementary abelian, and there exists a subgroup $T$ of $P$ such that $P=N\times T$. Then $T\cong P^@$, and $\Phi(P)=\Phi(T)$. Moreover $s\big(\Phi(P)\big)=\Phi(Q)$ as $P$ is a $p$-group, and $s\big(Z(P)\big)\leq Z(Q)$ as $s$ is surjective. It follows that $s(N)$ is an elementary abelian central subgroup of $Q$, so $s(N)\leq \Phi(Q)$ since $Q$ is atoric, by Lemma~\ref{atoric}. Now $s(P)=Q=s(N)s(T)$, thus $Q=\Phi(Q)s(T)$, and $s(T)=Q$, as was to be shown.\endpf
\begin{mth}{Notation} When $P$ is a finite $p$-group, and $N$ is a maximal normal subgroup of $P$ such that $N\cap\Phi(P)=\un$, we set $P^{@}=P/N$.
\end{mth}
By Proposition~\ref{Pzero}, the group $P^{@}$ does not depend on the choice of $N$, up to isomorphism: it is the greatest atoric quotient of $P$, in the sense that any atoric quotient of $P$ is isomorphic to a quotient of $P^{@}$. In particular $P^@$ is trivial if and only if $P$ is elementary abelian.
\begin{mth}{Proposition} \label{PatisoQat}Let $s:P\twoheadrightarrow Q$ be a surjective group homomorphism. Then $P^@\cong Q^@$ if and only if $\Ker(s)\cap\Phi(P)=\un$.
\end{mth}
\pf Let $E$ be a maximal normal subgroup of $P$ such that $E\cap\Phi(P)=\un$, and $T$ be a subgroup of~$P$ such that $P=E\times T$. Then $E$ is elementary abelian, and $\Phi(P)=\Phi(T)$. Let $\pi:Q\to Q^@$ be the canonical projection. By definition, we have $T\cong P^@$, and by Proposition~\ref{Pzero}, we have $\pi\circ s(T)=\nolinebreak Q^@$. Hence  $Q^@$ is a quotient of $P^@$, and $P^@\cong Q^@$ if and only if the map $\pi\circ s$ induces an isomorphism from $T$ to $Q^@$, that is if $\Ker(\pi\circ s)\cap T=\un$. This implies $\Ker(s)\cap T=\un$, hence $\Ker(s)\cap\Phi(P)=\un$. 
\par 
Conversely, if $\Ker(s)\cap\Phi(P)=\un$, then $\Ker(s)\cap\Phi(T)=\un$. Now the group $M=\Ker(s)\cap T$ is a normal subgroup of $T$ such that $M\cap \Phi(T)=\un$. Since $T$ is atoric, it follows from Lemma~\ref{atoric} that $M=\un$, hence $s(T)\cong T$. Now $Q=s(E)s(T)$, and $s(E)$ is a central elementary abelian subgroup of~$Q$, since $s$ is surjective. Let $F$ be a complement of $G=s(E)\cap s(T)$ in $s(E)$. Then $Q=(F\cdot G)s(T)=F\cdot s(T)$, thus $Q=F\times s(T)$ since $F$ is central in~$Q$. It follows that $s(T)$ is a quotient of $Q$. Since $s(T)\cong T\cong P^@$ is atoric, the group $P^@$ is isomorphic to a quotient of $Q^@$, thus $P^@\cong Q^@$.\endpf
\begin{mth}{Proposition} \label{sousquotient}Let $P$ be a finite $p$-group, and $Q$ be a subquotient of~$P$. Then $Q^@$ is a subquotient of $P^@$.
\end{mth}
\pf Let $(V,U)$ be a section of $P$ such that $V/U\cong Q$. Then $Q^@$ is isomorphic to a quotient of $V^@$, by Lemma~\ref{Pzero}. Hence it suffices to prove that $V^@$ is a subquotient of $P^@$.\par
Let $E$ be a maximal normal subgroup of $P$ such that $E\cap \Phi(P)=\un$, and $T$ be a subgroup of $P$ such that $P=E\times T$. Then $V\leq E\times T$, so there exist a subgroup $F$ of $E$, a subgroup $X$ of $T$, a group $Y$, and surjective group homomorphisms $\alpha:F\to Y$ and $\beta:X\to Y$ such that
$$V=\{(f,x)\in F\times X\mid \alpha(f)=\beta(x)\}\mpoint$$
Now $F\leq E$ is elementary abelian. If $(f,x), (f',x')\in V$, then $[(f,x),(f',x')]=(1,[x,x'])$, so $[V,V]\leq \un\times [X,X]$. Conversely if $x,x'\in X$, then there exist $f,f'\in F$ such that $\alpha(f)=\beta(x)$ and $\alpha(f')=\beta(x')$, i.e. $(f,x), (f',x')\in V$. Then $[(f,x),(f',x')]=(1,[x,x'])$, and it follows that $[V,V]=\un\times [X,X]$. Similarly, if $(f,x)\in V$, then $(f,x)^p=(1,x^p)$. Conversely, if $x\in X$, then there exists $f\in F$ such that $\alpha(f)=\beta(x)$, i.e. $(f,x)\in V$, and $(1,x^p)=(f,x)^p$. It follows that $\Phi(V)=\un\times\Phi(X)$.\par
Now $N=\Ker(\alpha)\times \un$ is a normal subgroup of $V$, and $N\cap\Phi(V)=\un$. By Proposition~\ref{PatisoQat}, it follows that $V^@\cong(V/N)^@$. Moreover the group homomorphism $(f,x)\in V\mapsto x\in X$ is surjective with kernel $N$, hence $V/N\cong X$. It follows that $V^@\cong X^@$ is a isomorphic to a quotient of the subgroup $X$ of $T\cong P^@$. Hence $V^@$ is a subquotient of $P^@$, as was to be shown.\endpf
\pagebreak[3]
\begin{mth}{Proposition}\label{ExL} Let $P$ be a finite $p$-group, let $N$ be a normal subgroup of $P$ such that $P/N\cong P^@$, and let $Q$ be a subgroup of $P$. The following are equivalent:
\begin{enumerate}
\item $Q^@\cong P^@$.
\item $QN=P$.
\item There exists a central elementary abelian subgroup $E$ of $P$ such that $P=EQ$.
\item There exists an elementary abelian subgroup $E$ of $P$ such that $P=E\times Q$.
\end{enumerate}
\end{mth}
\pf \framebox{$1\Rightarrow 2$} Suppose $Q^@\cong P^@$. We have $N\cap \Phi(T)=\un$, by Proposition~\ref{PatisoQat}. Moreover $\Phi(Q)\leq\Phi(P)$, as $P$ is a $p$-group. Setting $M=N\cap Q$, we have $M\cap\Phi(Q)=\un$, so $(Q/M)^@\cong Q^@\cong P^@$. But $\sur{Q}=Q/M\cong QN/N$ is a subgroup of $P/N\cong P^@$, and moreover there exists an elementary abelian subgroup $E$ of $\sur{Q}$ such that $\sur{Q}\cong E\times \sur{Q}^@\cong E\times P^@$. Hence $E=\un$ and $\sur{Q}\cong QN/N\cong P/N$, so $QN=P$, as was to be shown.\mpn
\framebox{$2\Rightarrow 3$} We have $N\cap \Phi(P)=\un$, by Proposition~\ref{PatisoQat}. Hence $N$ is elementary abelian, and central in $P$, and 2 implies 3.\mpn
\framebox{$2\Rightarrow 3$} Let $E$ be an elementary abelian central subgroup of $P$ such that $P=EQ$. Let $F$ be a complement of $E\cap Q$ in $E$. Then $F$ is elementary abelian and central in $P$. Moreover $QF=QE=P$, and $Q\cap F=\un$. Hence $P=F\times Q$.\mpn
\framebox{$4\Rightarrow 1$} If $P=E\times Q$ and $E$ is elementary abelian, then $\Phi(P)=\un\times \Phi(Q)$. Thus $E\cap \Phi(P)=\un$, so $(P/E)^@\cong P^@$ by Proposition~\ref{PatisoQat}, and $Q^@\cong P^@$.\endpf
\begin{mth}{Proposition} \begin{enumerate}
\item Let $L$ be an atoric $p$-group, let $P=E\times L$ and $Q=F\times L$, where $E$ and~$F$ are elementary abelian $p$-groups, and let $s:P\to Q$ be a group homomorphism. Then $s$ is surjective if and only if there exist a surjective group homomorphism $a:E\to F$, group homomorphisms $b:L\to F$ and $c:E\to\Omega_1Z(L)$, and an automorphism $d$ of $L$ such that
$$\forall (e,l)\in E\times L,\;\;\;s(e,l)=\big(a(e)b(l),c(e)d(l)\big)\mpoint$$
Moreover in this case $b\circ c(e)=1$ for any $e\in E$, and $s$ is an isomorphism if and only if $a$ is an isomorphism.
\item Let $P$ be a finite $p$-group. For a group homomorphism 
$$\lambda:P\to \Omega_1Z(P)\cap\Phi(P)\mvirg$$
let $\alpha_\lambda:P\to P$ be defined by $\alpha(x)=x\lambda(x)$, for $x\in P$. Then $\alpha_\lambda$ is an automorphism of $P$.
\item Let $P$ be a finite $p$-group, and let $P=E\times Q$, where $Q$ is atoric  and $E$ is elementary abelian. Then the correspondence $\lambda\mapsto \alpha_\lambda(E)$ is a bijection from the set of group homomorphisms $\lambda:P\to\Omega_1Z(P)\cap\Phi(P)$ such that $Q\leq \Ker\,\lambda$ to the set of subgroups $N$ of $P$ such that $P=N\times Q$. 
\end{enumerate}
\end{mth}
\pf (1) If $s$ is surjective, then $s(E)$ is central in $Q$, so $s(E)\leq\Omega_1Z(Q)=F\times\Omega_1Z(L)$. Hence there exists group homomorphisms $a:E\to F$ and $c:E\to\Omega_1Z(L)$ such that $s(e,1)=\big(a(e),c(e)\big)$, for any $e\in E$. Let $b: L\to F$ and $d:L\to L$ be the group homomorphisms defined by $s(1,l)=\big(b(l),d(l)\big)$, for $l\in L$. Then $s(e,l)=s(e,1)s(1,l)=\big(a(e)b(l),c(e)d(l)\big)$ for all $(e,l)\in P$. Moreover $b\circ c(e)=1$ for any $e\in E$, since $c(E)\leq\Omega_1Z(L)\leq\Phi(L)$, as $L$ is atoric, and $\Phi(L)\leq\Ker\,b$, as $F$ is elementary abelian.\par
Now the composition of $s$ with the projection $F\times L\to L$ is surjective, hence $s(\un\times L)=L$ by Proposition~\ref{Pzero}. In other words $d$ is surjective, hence it is an automorphism of $L$.\par
Since $s$ is surjective, for any $(f,y)\in Q$, there exists $(e,x)\in P$ such that $a(e)b(x)=f$ and $c(e)d(x)=y$. The latter gives $x=d^{-1}\big(c(e)^{-1}y\big)$. Then $b(x)=bd^{-1}\big(c(e)^{-1})bd^{-1}(y)$, and $bd^{-1}\big(c(e)^{-1}\big)=1$ since $d^{-1}\big(c(e)^{-1})\in d^{-1}\Omega_1Z(L)=\Omega_1Z(L)$, and $\Omega_1Z(L)\leq\Phi(L)\leq \Ker\,b$. Then $b(x)=bd^{-1}(y)$, and $f=a(e)bd^{-1}(y)$. In particular, taking $y=1$, we get that for any $f\in L$, there exists $e\in E$ such that $f=a(e)$. In other words $a$ is surjective.\par
Conversely, given a surjective group homomorphism $a:E\to F$, a group homomorphism $b:L\to F$, a group homomorphism $c:E\to\Omega_1Z(L)$, and an automorphism $d$ of $L$, we can define $s:P\to Q$ by $s(e,x)=\big(a(e)b(x),c(e)d(x)\big)$, for $(e,x)\in P$. This is clearly a group homomorphism, as $F$ is abelian, and the image of $c$ is central in $L$. We have again $\Omega_1Z(L)\leq\Phi(L)\leq\Ker\,b$, since $F$ is elementary abelian. If $(f,y)\in Q$, we can choose an element $e\in E$ such that $f=a(e)bd^{-1}(y)$, and then set $x=d^{-1}\big(c(e)^{-1}y\big)$, i.e. $c(e)d(x)=y$. We also have $b(x)=bd^{-1}(y)$, since $d^{-1}\big(c(e)\big)\in\Omega_1Z(L)$, so $f=a(e)b(x)$. Hence $s(e,x)=(f,y)$, and $s$ is surjective.\par
Finally if $s$ is an isomorphism, then $E\cong F$, and then the surjection $a$ is an isomorphism. Conversely, if $a$ is an isomorphism, then $E\cong F$, so $P\cong Q$, and the surjection $s$ is an isomorphism.\mpn
(2) Clearly $\alpha_\lambda$ is a group homomorphism, since $\lambda(P)\leq Z(P)$. Moreover if $x\in \Ker\,\alpha_\lambda$, then $\lambda(x)=x$, so $x\in\Omega_1Z(P)\cap\Phi(P)\leq \Phi(P)\leq\Ker\,\lambda$, since $\Omega_1Z(P)\cap\Phi(P)$ is elementary abelian. Thus $x=1$, and $\alpha_\lambda$ is injective. Hence it is an automorphism.\mpn
(3) Since $P=E\times Q$, we have $\Omega_1Z(P)=E\times \Omega_1Z(Q)$, and $\Phi(P)=\un\times\nolinebreak\Phi(Q)$. So if $\lambda$ is a group homomorphism from $P$ to $\Omega_1Z(P)\cap\Phi(P)$ with $Q\leq\Ker\,\lambda$, we have $\lambda(e,l)=\big(1,\beta(e)\big)$ for some group homomorphism $\beta:E\to\Omega_1Z(Q)$. Then the group $N=\alpha_\lambda(E)=\{\big(e,\beta(e)\big)\mid e\in E\}$ is central in $P$. Moreover $N\cap Q=\un$, and $NQ=P$, so $P=N\times Q$. Note that $N$ determines the homomorphism $\beta$, hence also the homomorphism $\lambda$, so the map $\lambda\mapsto \alpha_\lambda(E)$ is injective.\par
It is moreover surjective: indeed, if $N$ is a subgroup of $P=E\times Q$ such that $P=N\times Q$, then $N\cong P/Q\cong E$ is elementary abelian, hence central in~$P$. Since $NQ=P$, for any $e\in E$, there exists $(a,b)\in N$ and $q\in Q$ such that $(e,1)=(a,b)(1,q)$, that is $e=a$ and $q=b^{-1}$. In other words $p_1(N)=E$. Moreover $N\cap Q=\un$, so $k_2(N)=\un$. So for $e\in E$, there exists a unique $x\in Q$ such that $(e,x)\in N$. Setting $x=\beta(e)$, we get a group homomorphism $\beta:E\to Q$, such that $N=\big\{\big(e,\beta(e)\big)\mid e\in E\big\}$. Since $N$ is central in $P$, the image of $\beta$ is contained in $\Omega_1Z(Q)\leq\Phi(Q)$. Moreover $\Omega_1Z(P)=E\times \Omega_1Z(Q)$, and $\Phi(P)=1\times\Phi(Q)$, so $\big(\un\times\beta(E)\big)\leq\Omega_1Z(P)\cap\Phi(P)$. Setting $\lambda(e,l)=\big(1,\beta(e)\big)$, we get a group homomorphism from $P$ to $\Omega_1Z(P)\cap \Phi(P)$, such that $Q\leq \Ker\,\lambda$, and $N=\alpha_\lambda(E)$.\endpf
\section{Splitting the biset category of $p$-groups, when $p\in\Cf^\times$}
\begin{mth}{Notation and Definition} Let $\Cf\mathcal{C}_p$ denote the full subcategory of the biset category $\Cf\mathcal{C}$ consisting of finite $p$-groups. A {\em $p$-biset functor} over~$\Cf$ is an $\Cf$-linear functor from $\Cf\mathcal{C}_p$ to the category of $\Cf$-modules. Let $\CF_{p,\Cf}$ denote the full subcategory of $\CF_\Cf$ consisting of $p$-biset functors over $\Cf$.
\end{mth}
In the statements below, we indicate by {\rm [ }$p\in \Cf^\times${\rm ] } the assumption that $p$ is invertible in $R$.
\begin{mth}{Theorem} \label{epsilon}{\rm [ }$p\in \Cf^\times${\rm ] } Let $P$ and $Q$ be finite $p$-groups, let $(T,S)$ be a minimal section of $P$, and $(V,U)$ be a minimal section of $Q$. Then
$$\epsilon_{V,U}^Q\,\Cf B(Q,P)\,\epsilon_{T,S}^P\neq\zero \;\implies\;(V/U)^@\cong (T/S)^@\mpoint$$
\end{mth}
\pf If $\epsilon_{V,U}^Q\Cf B(Q,P)\epsilon_{T,S}^P$, there exists $a\in \Cf B(Q,P)$ such that 
$$\epsilon_{V,U}^Q\,a\,\epsilon_{T,S}^P=\Indinf_{V/U}^Q\varphi_\un^{V/U}\Defres_{V/U}^Q\,a\,\Indinf_{T/S}^P\varphi_{1}^{T/S}\Defres_{T/S}^P\neq 0\mvirg$$
and in particular the element $b=\Defres_{V/U}^Q\,a\,\Indinf_{T/S}^P$ of $\Cf B(V/U,T/S)$ is such that $\varphi_\un^{V/U}\,b\,\varphi_\un^{T/S}\neq 0$. It follows that there is a subgroup $L$ of the product $(V/U)\times(T/S)$ such that
$$\varphi_\un^{V/U}\big[\big((V/U)\times(T/S)\big)/L\big]\varphi_\un^{T/S}\neq 0\mpoint$$
Then Theorem~\ref{idempotents E} implies that $p_1(L)=V/U$, $k_1(L)\cap \Phi(V/U)=\un$, $p_2(L)=T/S$, and $k_2(L)\cap \Phi(T/S)=\un$. By Proposition~\ref{PatisoQat}, it follows that
$$(V/U)^@\cong \big(p_1(L)/k_1(L)\big)^@\cong \big(p_2(L)/k_2(L)\big)^@\cong(T/S)^@\mvirg$$
as was to be shown.\endpf
\begin{mth}{Notation}\label{bLP} {\rm [ }$p\in \Cf^\times${\rm ] } Let $L$ be an atoric $p$-group. If $P$ is a finite $p$-group, we set
$$b_L^P=\sumb{(T,S)\in[\mathcal{M}(G)]}{\rule{0ex}{1.6ex}(T/S)^@\cong L}\epsilon_{T,S}^P\mpoint$$
\end{mth}
\begin{mth}{Theorem} \label{b_L}{\rm [ }$p\in \Cf^\times${\rm ] } 
\begin{enumerate}
\item Let $L$ be an atoric $p$-group, and $P$ be a finite $p$-group. Then $b_L^P\neq 0$ if and only if $L\sqsubseteq P^@$.
\item Let $L$ and $M$ be atoric $p$-groups, and let $P$ and $Q$ be finite $p$-groups. If $b_M^Q\Cf B(Q,P)b_L^P\neq\zero$, then $M\cong L$.
\item Let $L$ be an atoric $p$-group, and let $P$ and $Q$ be finite $p$-groups. Then for any $a\in \Cf B(Q,P)$
$$b_L^Q\,a=a\,b_L^P\mpoint$$
\item The family of elements $b_L^P\in\Cf B(P,P)$, for finite $p$-groups $P$, is an idempotent endomorphism $b_L$ of the identity functor of the category $\Cf\CC_p$ (i.e. an idempotent of the center of $\Cf\CC_p$). The idempotents $b_L$, for $L\in[\At_p]$, are orthogonal, and their sum is equal to the identity element of the center of $\Cf\CC_p$.
\item For a given finite $p$-group $P$, the elements $b_L^P$, for $L\in [\At_p]$ such that $L\sqsubseteq P^@$, are non zero orthogonal central idempotents of $\Cf B(P,P)$, and their sum is equal to the identity of $\Cf B(P,P)$.
\end{enumerate}
\end{mth}
\pf (1) The idempotent $b_L^P$ is non zero if and only if there exists a minimal section $(T,S)$ of $P$ such that $(T/S)^@\cong L$. Then $L\sqsubseteq P^@$, by Proposition~\ref{sousquotient}. Conversely, if $L\sqsubseteq P^@$, then $L\sqsubseteq P$, and there exists a minimal section $(T,S)$ of $P$ such that $T/S\cong L$. Then $(T/S)^@\cong L^@\cong L$, so $\epsilon_{T,S}^P$ appears in the sum defining $b_L^p$, thus $b_L^P\neq 0$.\mpn
(2) If $b_M^Q\Cf B(Q,P)b_L^P\neq\zero$, then there exist a minimal section $(V,U)$ of $Q$ with $(V/U)^@\cong M$ and a minimal section $(T,S)$ of $P$ with $(T/S)^@\cong L$ such that $\epsilon_{V,U}^Q\Cf B(Q,P)\epsilon_{T,S}^P\neq 0$. Then $(V/U)^@\cong (T/S)^@$ by Theorem~\ref{epsilon}, that is $M\cong L$.\mpn
(3) The identity element of $\Cf B(P,P)$ is equal to the sum of the idempotents $\epsilon_{T,S}^P$, for $(T,S)\in[\CM(P)]$. Grouping those idempotents $\epsilon_{T,S}^P$ for which $(T/S)^@$ is isomorphic to a given $L\in[\At_p]$ shows that the identity element of $\Cf B(P,P)$ is equal to the sum of the idempotents $b_L^P$, for $L\in[\At_p]$ (and there are finitely many non zero $b_L^P$, by (1)). It follows that
\begin{eqnarray*}
b_M^Q\,a&=&b_M^Q\,a\,\sum_{L\in [\At_p]}b_L^P=\sum_{L\in[\At_p]}b_M^Q\,a\,b_L^P\\
&=&b_M^Q\,a\,b_M^P \;\;\hbox{[by (2)]}\\
&=&\sum_{L\in[\At_p]}b_L^Q\,a\,b_M^P\;\;\hbox{[by (2)]}\\
&=&a\,b_{M}^P\mvirg
\end{eqnarray*}
since $\sum_{L\in[\At_p]}\limits b_L^Q$ is the identity element of $\Cf B(Q,Q)$.\par
It follows that the family $b_L^P$, where $P$ is a finite $p$-group, is an element $b_L$ of the center of $\Cf\CC_p$. Clearly $b_L^2=b_L$, and if $L$ and $M$ are non isomorphic atoric $p$-groups, then $b_Lb_M=0$, by (2). Moreover the infinite sum $\sum_{L\in[\At_p]}\limits b_L$ is actually locally finite, i.e. for each finite $p$-group $P$, the sum $\sum_{L\in[\At_p]}\limits b_L^P$ has only finitely many non zero terms. The sum $\sum_{L\in[\At_p]}\limits b_L$ is clearly equal to the identity endomorphism of the identity functor of $\Cf\CC_p$. \mpn
(4) This is a straightforward consequence of (1) and (3).\endpf
\begin{mth}{Corollary}\label{decomposition}{\rm [ }$p\in \Cf^\times${\rm ] } 
\begin{enumerate} 
\item Let $L$ be an atoric $p$-group. For a $p$-biset functor $F$, the family of maps $F(b_L^P):F(P)\to F(P)$, for finite $p$-groups $P$, is an endomorphism of~$F$, denoted by $F(b_L)$.
\item If $\theta:F\to G$ is a natural transformation of $p$-biset functors, the diagram
$$\xymatrix{
F\ar[d]_-\theta\ar[r]^-{F(b_L)}&F\ar[d]^-\theta\\
G\ar[r]_-{G(b_L)}&G
}
$$
is commutative. Hence the family of endomorphisms $F(b_L)$, for $p$-biset functors $F$, is an idempotent of the center of the category $\CF_{p,\Cf}$, denoted by $\widehat{b}_L$.
\item The idempotents $\widehat{b}_L$, for $L\in[\At_p]$, are orthogonal idempotents of the center of $\CF_{p,\Cf}$, and their sum is the identity. 
\item If $F$ is a $p$-biset functor over $\Cf$, let $\widehat{b}_LF$ denote the image of the endomorphism $F(b_L)$ of $F$. Then $F=\dirsum{L\in[\At_p]}\widehat{b}_LF$.
\item Let $\widehat{b}_L\CF_{p,\Cf}$ denote the full subcategory of $\CF_{p,\Cf}$ consisting of functors $F$ such that $F=\widehat{b}_LF$. Then $\widehat{b}_L\CF_{p,\Cf}$ is an abelian subcategory of $\CF_{p,\Cf}$. Moreover the functor
$$F\in\CF_{p,\Cf}\mapsto (\widehat{b}_LF)_{L\in[\At_p]}\in\prod_{L\in[\At_p]}\widehat{b}_L\CF_{p,\Cf}$$
is an equivalence of categories.
\end{enumerate}
\end{mth}
\pf All assertions are straightforward consequences of Theorem~\ref{b_L}.\endpf
\begin{mth}{Notation} For an atoric $p$-group $L$, let $\Cf\CC_p^L$ denote the full subcategory of $\Cf\CC_p$ consisting of the class $\CY_L$ of finite $p$-groups $P$ such that $P^@\sqsubseteq L$. When $p\in R^\times$, Let moreover
$$b_L^+=\sumb{H\in[\At_p]}{H\sqsubseteq L}b_H$$
be the sum of the idempotents $b_H$ corresponding to atoric subquotients of $L$, up to isomorphism.
\end{mth}
The class $\CY_L$ is closed under taking subquotients, by Proposition~\ref{sousquotient}. It follows that we can apply the results of Section 6 (Appendix) of~\cite{both3}: if $F$ is a $p$-biset functor over $\Cf$, we can restrict $F$ to an $\Cf$-linear functor from $\Cf\CC_p^L$ to $\gMod{\Cf}$. This yields a forgetful functor $\CO_{\CY_L}: \CF_{p,\Cf}\to \mathsf{Fun}_\Cf\big(\Cf\CC_p^L,\gMod{\Cf})$. The right adjoint $\CR_{\CY_L}$ of this functor is described in full detail in Section~6 of~\cite{both3}, as follows: if $G$ is an $\Cf$-linear functor from $\Cf\CC_p^L$ to $\gMod{\Cf}$, and $P$ is a finite $p$-group, set
\begin{equation}\label{inverse limit}
\CR_{\CY_L}(G)(P)=\limproj{(X,M)\in\Sigma_L(P)}G(X/M)
\end{equation}
the inverse limit of modules $G(X/M)$ on the set $\Sigma_L(P)$ of sections $(X,M)$ of~$P$ such that $(X/M)^@\sqsubseteq L$, i.e. the set of sequences $(l_{X,M})_{(X,M)\in \Sigma_{L}(P)}$ with the following properties:
\begin{enumerate}
\item if $(X,M)\in\Sigma_L(P)$, then $l_{X,M}\in G(X/M)$.
\item if $(X,M), (Y,N)\in \Sigma_L(P)$ and $M\leq N\leq Y\leq X$, then 
$$\Defres_{Y/N}^{X/M}l_{X,M}=l_{Y,N}\mpoint$$
\item if $x\in P$ and $(X,M)\in\Sigma_L(P)$, then ${^xl}_{X,M}=l_{^xX,^xM}$.
\end{enumerate}
Recall now that for finite groups $P$ and $Q$, and for a finite $(Q,P)$-biset $U$, for a subgroup $T$ of $Q$ and an element $u$ of $U$, the subgroup $T^u$ of $P$ is defined by $T^u=\{x\in P\mid\exists t\in T\;\;tu=ux\}$. By Lemma~6.4 of~\cite{both3}, if $(T,S)$ is a section of $Q$, then $(T^u,S^u)$ is a section of $P$, and $T^u/S^u$ is a subquotient of $T/S$. \par
With this notation, when $P$ and $Q$ are finite $p$-groups, when $U$ is a finite $(Q,P)$-biset, and $l= (l_{X,M})_{(X,M)\in\Sigma_L(P)}$ is an element of $\CR_{\CY_L}(G)(P)$, we denote by $Ul$ the sequence indexed by $\Sigma_L(Q)$ defined by
$$(Ul)_{Y,N}=\sum_{u\in [Y\dom U/P]_L}(N\dom Yu)(l_{Y^u,N^u})$$
where $[Y\dom U/P]$ is a set of representatives of $(Y\times P)$-orbits on $U$, and $N\dom Yu$ is viewed as a $(Y/N,Y^u/N^u)$-biset. It shown in Section~6 of~\cite{both3} that $Ul\in \CR_{\CY_L}(G)(Q)$, and that $\CR_{\CY_L}(G)$ becomes a $p$-biset functor in this way. Moreover\footnote{In Theorem~6.15 of \cite{both3}, only the case $\Cf=\Z$ is considered, but the proofs extend trivially to the case of an arbitrary commutative ring $\Cf$}:
\pagebreak[3]
\begin{mth}{Theorem}{\rm [\cite{both3} Theorem~6.15]}\label{appendix} The assignment $G\mapsto \CR_{\CY_L}(G)$ is an $\Cf$-linear functor $\CR_{\CY_L}$ from $\mathsf{Fun}_\Cf\big(\Cf\CC_p^L,\gMod{\Cf}\big)$ to $\CF_{p,\Cf}$, which is right adjoint to the forgetful functor $\CO_{\CY_L}$. Moreover the composition $\CO_{\CY_L}\circ \CR_{\CY_L}$ is isomorphic to the identity functor of $\mathsf{Fun}_\Cf\big(\Cf\CC_p^L,\gMod{\Cf}\big)$.
\end{mth}
\pagebreak[3]
\begin{mth}{Theorem}\label{equivalence}{\rm [ }$p\in \Cf^\times${\rm ] }  For an atoric $p$-group $L$, let $\widehat{b}_L^+\CF_{p,\Cf}$ be the full subcategory of $\CF_{p,\Cf}$ consisting of functors $F$ such that $\widehat{b}_L^+F=F$. Then the forgetful functor $\CO_{\CY_L}$ and its right adjoint $\CR_{\CY_L}$ restrict to quasi-inverse equivalences of categories
$$\xymatrix{
\widehat{b}_L^+\CF_{p,\Cf}\ar[r]<.5ex>^-{\CO_{\CY_L}}&\ar[l]<.5ex>^-{\CR_{\CY_L}} \mathsf{Fun}_\Cf\big(\Cf\CC_p^L,\gMod{\Cf}\big)\mpoint
}$$
\end{mth}
\pf \sou{First step:} The first thing to check is that the image of the functor $\CR_{\CY_L}$ is contained in $\widehat{b}_L^+\CF_{p,\Cf}$. We first prove that if $H$ is an atoric $p$-group, if $F\in\CF_{p,\Cf}$, and if $\CO_{\CY_L}(\widehat{b}_HF)\neq 0$, then $H\sqsubseteq L$: indeed in that case, there exists $P\in\CY_L$ such that $b_H^PF(P)\neq 0$. In particular $b_H^P\neq 0$ by Theorem~\ref{b_L}, hence $H\sqsubseteq P^@$. Since $P^@\sqsubseteq L$ as $P\in\CY_L$, it follows that $H\sqsubseteq L$, as claimed.\par
In particular 
$$\CO_{\CY_L}(F)=\CO_{\CY_L}\big(\sumb{H\in[\At_p]}{H\sqsubseteq L}\widehat{b}_HF\big)=\CO_{\CY_L}\big(\widehat{b}_L^+F\big)\mpoint$$
Set $\CG_p^L=\mathsf{Fun}_\Cf\big(\Cf\CC_p^L,\gMod{\Cf}\big)$, and let $G\in\CG_p^L$. Let $H$ be an atoric $p$-group such that $H\not\sqsubseteq L$. If $F\in\CF_{p,\Cf}$, then
\begin{eqnarray*}
\Hom_{\CF_{p,\Cf}}\big(F,\widehat{b}_H\CR_{\CY_L}(G)\big)&=&\Hom_{\CF_{p,\Cf}}\big(\widehat{b}_HF,\widehat{b}_H\CR_{\CY_L}(G)\big)\\
&=&\Hom_{\CF_{p,\Cf}}\big(\widehat{b}_HF,\CR_{\CY_L}(G)\big)\\
&\cong&\Hom_{\CG_p^L}\big(\CO_{\CC_{\CY_L}}\big(\widehat{b}_HF\big),G\big)=\zero\mpoint\\
\end{eqnarray*}
So the functor $F\mapsto\Hom_{\CF_{p,\Cf}}\big(F,\widehat{b}_H\CR_{\CY_L}(G)\big)$ is the zero functor, and it follows from Yoneda's lemma that $\widehat{b}_H\CR_{\CY_L}(G)=0$ if $H\not\sqsubseteq L$. In other words $\CR_{\CY_L}(G)=\widehat{b}_L^+\CR_{\CY_L}(G)$, as was to be shown.\mpn
\sou{Second step:} The first step shows that we have adjoint functors
$$\xymatrix{
\widehat{b}_L^+\CF_{p,\Cf}\ar[r]<.5ex>^-{\CO_{\CY_L}}&\ar[l]<.5ex>^-{\CR_{\CY_L}} \mathsf{Fun}_\Cf\big(\Cf\CC_p^L,\gMod{\Cf}\big)=\CG_p^L\mpoint
}$$
Moreover, the composition $\CO_{\CY_L}\circ\CR_{\CY_L}$ is isomorphic to the identity functor, by Theorem~\ref{appendix}. All we have to show is that the unit of the adjunction is also an isomorphism, in other words, that for any $F\in\widehat{b}_L^+\CF_{p,\Cf}$ and any finite $p$-group $P$, the natural map
\begin{equation}\label{unite}
\eta_P:F(P)\to\CR_{\CY_L}\CO_{\CY_L}(F)(P)=\limproj{(X,M)\in\Sigma_L(P)}F(X/M)
\end{equation}
sending $u\in F(P)$ to the sequence $\big(\Defres_{X/M}^Pu\big)_{(X,M)\in\Sigma_L(P)}$, is an isomorphism. \par
The map $\eta_P$ is injective: indeed, if $u\in F(P)$, then $u=\sumb{H\in[\At_p]}{H\subseteq L}\limits b_H^Pu$, as $F=\widehat{b}_L^+F$. If $\Defres_{X/M}^Pu=0$ for any section $(X,M)$ of $P$ with $(X/M)^@\sqsubseteq L$, then $F(\epsilon_{T,S}^P)(u)=0$ for any section $(T,S)$ of $P$ such that $(T/S)^@\sqsubseteq L$, by Proposition~\ref{formule epsilon} and Proposition~\ref{sousquotient}. In particular $b_H^Pu=0$ for any atoric subquotient $H$ of $L$, hence $u=0$.\par
To prove that $\eta_P$ is also surjective, we generalize the construction of Theorem A.2 of~\cite{both2} (which is the case $L=\un$), and we define, for an element $v=(v_{X,M})_{(X,M)\in\Sigma_L(P)}$ in $\CR_{\CY_L}\CO_{\CY_L}(F)(P)$, an element $u=\iota_P(v)$ of $F(P)$ by 
$$u=\frac{1}{|P|}\sumb{(T,S)\in\CM(P)}{(T/S)^@\sqsubseteq L}\sumb{X\leq T,\, M\normal T}{S\leq M\leq\Phi(T)\leq X\leq T}|X|\mu(X,T)\mu_{\normal T}(S,M)\Indinf_{X/M}^Pv_{X,M}\mpoint$$
This yields an $\Cf$-linear map $\iota_P:\CR_{\CY_L}\CO_{\CY_L}(F)(P)\to F(P)$. \par
For $(Y,N)\in\Sigma_L(P)$, set $u_{Y,N}=\Defres_{Y/N}^Pu$. Then:
$$u_{Y,N}=\!\!\!\sumb{(T,S)\in\CM(P)}{(T/S)^@\sqsubseteq L}\!\sumb{X\leq T,\, M\normal T}{S\leq M\leq\Phi(T)\leq X\leq T}\!\!\!\!\!\!\frac{|X|}{|P|}\mu(X,T)\mu_{\normal T}(S,M)\Defres_{Y/N}^P\Indinf_{X/M}^Pv_{X,M}.$$
Moreover
$$\Defres_{Y/N}^P\Indinf_{X/M}^Pv_{X,M}=\sum_{g\in[Y\dom P/X]}\Indinf_{J_g/J'_g}^{Y/N}\Iso(\phi_g)\Defres_{I_g/I'_g}^{^gX/^gM}{^gv}_{X,N}\mvirg$$
where $J_g=N(Y\cap {^gX})$, $J'_g=N(Y\cap {^gM})$, $I_g={^gM}(Y\cap {^gX})$, $I'_g={^gM}(N\cap {^gX})$, and $\phi_g$ is the isomorphism $I_g/I'_g\to J_g/J'_g$ sending $xI'_g$ to $xJ'_g$, for $x\in Y\cap {^gX}$. Hence
\begin{eqnarray*}
\Defres_{Y/N}^P\Indinf_{X/M}^Pv_{X,M}&=&\sum_{g\in[Y\dom P/X]}\Indinf_{J_g/J'_g}^{Y/N}\Iso(\phi_g)v_{I_g,I'_g}\\
&=&\frac{|Y\cap{^gX}|}{|Y||X|}\sum_{g\in P}\Indinf_{J_g/J'_g}^{Y/N}\Iso(\phi_g)v_{I_g,I'_g}\mpoint
\end{eqnarray*}
Thus
$$u_{Y,N}=\!\!\!\sum_{\substack{(T,S)\in\CM(P)\\(T/S)^@\sqsubseteq L\\X\leq T,\, M\normal T\\S\leq M\leq\Phi(T)\leq X\leq T\\g\in P}}\!\!\!\!\!\!\frac{|Y\cap{^gX}|}{|P||Y|}\mu(X,T)\mu_{\normal T}(S,M)\Indinf^{Y/N}_{J_g/J'_g}\Iso(\phi_g)v_{I_g,I'_g}.$$
Now $\mu(X,T)=\mu({^gX},{^gT})$ and $\mu_{\normal T}(S,M)=\mu_{\normal{^gT}}({^gS},{^gM})$, so summing over $({^gT},{^gS},{^gX},{^gM})$ instead of $(T,S,X,M)$ we get
$$u_{Y,N}=\!\!\!\sum_{\substack{(T,S)\in\CM(P)\\(T/S)^@\sqsubseteq L\\X\leq T,\, M\normal T\\S\leq M\leq\Phi(T)\leq X\leq T}}\!\!\!\!\!\!\frac{|Y\cap{X}|}{|Y|}\mu(X,T)\mu_{\normal T}(S,M)\Indinf^{Y/N}_{J_1/J'_1}\Iso(\phi_1)v_{I_1,I'_1}.$$
Setting $W=Y\cap X$, we have $J_1=NW$, $J'_1=N(W\cap M)$, $I_1=MW$, $I'_1=M(N\cap W)$, and these four groups only depend on $W$, once $M$ and $N$ are given. Hence, for given $T,S$ and $M$, we can group together the terms of the above summation for which $Y\cap X$ is a given subgroup $W$ of $Y\cap T$. This gives
$$u_{Y,N}=\!\!\!\sum_{\substack{(T,S)\in\CM(P)\\(T/S)^@\sqsubseteq L\\M\normal T\\S\leq M\leq\Phi(T)\\W\leq Y\cap T}}\big(\sum_{\substack{\Phi(T)\leq X\leq T\\X\cap Y=W}}\mu(X,T)\big)\frac{|W|}{|Y|}\mu_{\normal T}(S,M)\Indinf^{Y/N}_{J_1/J'_1}\Iso(\phi_1)v_{I_1,I'_1}.$$
Moreover $\sum_{\substack{\Phi(T)\leq X\leq T\\X\cap Y=W}}\limits\mu(X,T)=\sum_{\substack{X\leq T\\X\cap (Y\cap T)=W}}\limits\mu(X,T)$, since $\mu(X,T)=0$ unless $X\geq \Phi(T)$, and the latter summation vanishes unless $Y\cap T=T$, by classical combinatorial lemmas (\cite{stanley} Corollary 3.9.3). This gives:
$$u_{Y,N}=\!\!\!\sum_{\substack{(T,S)\in\CM(P)\\(T/S)^@\sqsubseteq L\\M\normal T\\S\leq M\leq\Phi(T)\leq W\leq T\leq Y}}\frac{|W|}{|Y|}\mu(W,T)\mu_{\normal T}(S,M)\Indinf^{Y/N}_{J_1/J'_1}\Iso(\phi_1)v_{I_1,I'_1}.$$
Moreover in this summation $J_1=NW$, $J'_1=N(W\cap M)=NM$, $I_1=MW=W$, $I'_1=M(N\cap W)=MN\cap W$. All these groups remain unchanged if we replace $M$ by $M\big(N\cap\Phi(T)\big)$, so for given $T,S$ and $W$, we can group together those terms for which $M\big(N\cap \Phi(T)\big)$ is a given normal subgroup $U$ of $T$ with $U\leq\Phi(T)$. The sum $\sum_{\substack{S\leq M\normal T\\M\big(N\cap\Phi(T)\big)=U}}\limits\mu_{\normal T}(S,M)$ is equal to 0 (by the same above-mentioned classical combinatorial lemmas) unless $N\cap\Phi(T)\leq S$. Hence
$$u_{Y,N}=\!\!\!\!\!\sum_{\substack{(T,S)\in\CM(P)\\(T/S)^@\sqsubseteq L\\U\normal T\\N\cap \Phi(T)\leq S\leq U\leq\Phi(T)\leq W\leq T\leq Y}}\!\!\!\!\frac{|W|}{|Y|}\mu(W,T)\mu_{\normal T}(S,U)\Indinf^{Y/N}_{J_1/J'_1}\Iso(\phi_1)v_{I_1,I'_1},$$
where $J_1=NW$, $J'_1=NU$, $I_1=W$, $I'_1=UN\cap W$.\par
Now if $N\cap \Phi(T)\leq S\leq\Phi(T)\leq T\leq Y$, then $(TN/N)^@\sqsubseteq (Y/N)^@$. Moreover the normal subgroup $(N\cap T)/\big(N\cap\Phi(T)\big)$ of $T/\big(N\cap\Phi(T)\big)$ intersects trivially the Frattini subgroup 
$$\Phi\Big(T/\big(N\cap\Phi(T)\big)\Big)=\Phi(T)\big(N\cap\Phi(T)\big)/\big(N\cap\Phi(T)\big)\mvirg$$
so $\Big(T/\big(N\cap\Phi(T)\big)\Big)^@\cong (T/(T\cap N)^@\cong (TN/N)^@$ by Proposition~\ref{PatisoQat}. Then $(T/S)^@\sqsubseteq \Big(T/\big(N\cap\Phi(T)\big)\Big)^@\sqsubseteq (TN/N)^@\sqsubseteq (Y/N)^@$. As $(Y/N)^@\sqsubseteq L$ by assumption, it follows that
$$u_{Y,N}=\!\!\!\!\!\sum_{\substack{S\normal T\leq Y\\U\normal T\\N\cap \Phi(T)\leq S\leq U\leq\Phi(T)\leq W\leq T\leq Y}}\!\!\!\!\frac{|W|}{|Y|}\mu(W,T)\mu_{\normal T}(S,U)\Indinf^{Y/N}_{J_1/J'_1}\Iso(\phi_1)v_{I_1,I'_1}.$$
Now the sum $\sum_{\substack{S\normal T\\N\cap \Phi(T)\leq S\leq U}}\limits\mu_{\normal T}(S,U)$ is equal to zero unless $U=N\cap\Phi(T)$. Hence
$$u_{Y,N}=\sum_{\Phi(T)\leq W\leq T\leq Y}\limits\frac{|W|}{|Y|}\mu(W,T)\Indinf^{Y/N}_{J_1/J'_1}\Iso(\phi_1)v_{I_1,I'_1}.$$
For a given subgroup $W$ of $Y$, the sum $\sum_{\Phi(T)\leq W\leq T\leq Y}\limits\mu(W,T)$ is equal to $\sum_{W\leq T\leq Y}\limits\mu(W,T)$ since $\mu(W,T)=0$ unless $W\geq\Phi(T)$, and the latter is equal to zero if $W\neq Y$, and to 1 if $W=Y$. Thus
$$u_{Y,N}=\frac{|Y|}{|Y|}\Indinf^{Y/N}_{J_1/J'_1}\Iso(\phi_1)v_{I_1,I'_1},$$
where $J_1=NY=Y$, $J'_1=N(Y\cap U)=N$, $I_1=Y=, I'_1=UN\cap Y=N$. Hence $I_1=J_1=Y$ and $I'_1=J'_1=N$, so $\phi_1$ is equal to the identity. It follows that $u_{Y,N}=v_{Y,N}$ for any $(Y,N)\in\Sigma_L(P)$, so $\eta_P(u)=v$. This proves that the map $\eta_P$ is surjective, hence an isomorphism, with inverse $\iota_P$. This completes the proof of Theorem~\ref{equivalence}.\endpf
\begin{mth}{Definition}\label{CCpdeL}
Let $\Cf\CC_p^{\sharp L}$ be the following category:
\begin{itemize}
\item The objects of $\Cf\CC_p^{\sharp L}$ are the finite $p$-groups $P$ such that $P^@\cong L$.
\item If $P$ and $Q$ are finite $p$-groups such that $P^@\cong Q^@\cong L$, then
$$\Hom_{\Cf\CC_p^{\sharp L}}(P,Q)=\Cf B(Q,P)/\sum_{L\not\sqsubseteq S}\Cf B(Q,S)B(S,P)$$
is the quotient of $\Cf B(Q,P)$ by the $\Cf$-submodule generated by all morphisms from $P$ to $Q$ in $\Cf\CC_p$ which factor through a $p$-group $S$ which do not admit $L$ as a subquotient.
\item The composition of morphisms in $\Cf\CC_p^{\sharp L}$ is induced by the composition of morphisms in $\Cf\CC_p$. 
\end{itemize}
\end{mth}
\begin{rem}{Remark} Morphisms in $\Cf\CC_p$ which factor through a $p$-group $S$ such that $L\not\sqsubseteq S$ clearly generate a two-sided ideal, so the composition in $\Cf\CC_p^{\sharp L}$ is well defined. Moreover the category $\Cf\CC_p^{\sharp L}$ is $\Cf$-linear. Let $\mathsf{Fun}_\Cf\big(\Cf\CC_p^{\sharp L},\gMod{\Cf}\big)$ denote the category of $\Cf$-linear functors from $\Cf\CC_p^{\sharp L}$ to the category $\gMod{\Cf}$ of $\Cf$-modules.
\end{rem}
\begin{mth}{Lemma}\label{abelemxL} Let $p$ be a prime, and $L$ be an atoric $p$-group. Let $P$ and $Q$ be finite $p$-groups.
\begin{enumerate}
\item If $P^@\cong L$ or $Q^@\cong L$, and if $M\leq (Q\times P)$, then $q(M)^@\sqsubseteq L$. Moreover $q(M)^@\cong L$ if and only if $L\sqsubseteq q(M)$.
\item If $P^@\cong Q^@\cong L$, then
$$\Hom_{\Cf\CC_p^{\sharp L}}(P,Q)=\Cf B(Q,P)/\sum_{S^@\sqsubset L}\Cf B(Q,S)B(S,P)$$
is also the quotient of $\Cf B(Q,P)$ by the $\Cf$-submodule generated by all morphisms from $P$ to $Q$ in $\Cf\CC_p$ which factor through a $p$-group $S$ such that $S^@$ is a proper subquotient of $L$.
\item If $P^@\cong Q^@\cong L$, then $\Hom_{\Cf\CC_p^{\sharp L}}(P,Q)$ has an $\Cf$-basis consisting of the (images of the) transitive $(Q,P)$-bisets $(Q\times P)/M$, where $M$ is a subgroup of $(Q\times P)$ such that $q(M)^@\cong L$ (up to conjugation).
\end{enumerate}
\end{mth}
\pf (1) Indeed $q(M)$ is a subquotient of $P$, and a subquotient of $Q$. Hence $q(M)^@$ is a subquotient of $P^@$ and a subquotient of $Q^@$, thus $q(M)\sqsubseteq L^@\cong L$. Now suppose that $q(M)^@\cong L$. Then $L$ is a quotient of $q(M)$, so $L\sqsubseteq q(M)$. Conversely, if $L\sqsubseteq q(M)$, then $L\cong L^@$ is a subquotient of $q(M)^@$, which is a subquotient of $L$. So $q(M)^@\cong L$.\mpn
(2) Let $S$ be a finite $p$-group such that $L\not\sqsubseteq S$, or equivalently $L\not\sqsubseteq S^@$. Any element of $\Cf B(Q,S)B(S,P)$ is a linear combination of $(Q,P)$-bisets of the form $(Q\times P)/(M*N)$, for $M\leq (Q\times S)$ and $N\leq (S\times P)$. This biset $(Q\times P)/(M*N)$ also factors though $T=q(M*N)$, by~\ref{factorize}. Moreover $T$ is a subquotient of $q(M)$ and $q(N)$, hence a subquotient of $Q$, $S$, and $P$. Hence $T^@\sqsubseteq Q^@\cong L$, and $T^@\ncong L$, since $L\not\sqsubseteq S^@$. Hence $T\sqsubset L$.\mpn
(3) The (images of the) elements $(Q\times P)/M$, where $M$ is a subgroup of $(Q\times P)$ such that $q(M)^@\cong L$ (up to conjugation), clearly generate $\Hom_{\Cf\CC_p^{\sharp L}}(P,Q)$. Moreover, the proof of (2) shows that they are linearly independent, since any transitive $(Q,P)$-biset $(Q\times P)/N$ appearing in an element of the sum $\sum_{S^@\sqsubset L}\limits\Cf B(Q,S)B(S,P)$ is such that $q(N)^@\sqsubset L$.
.\endpf
\begin{rem}{Remark} If $G$ is an $\Cf$-linear functor from $\Cf\CC_p^{\sharp L}$ to the category $\gMod{\Cf}$ of $\Cf$-modules, we can extend $G$ to an $\Cf$-linear functor from $\Cf\CC_p^L$ to $\gMod{\Cf}$ by setting $G(P)=\zero$ if $P$ is a finite $p$-group such that $P^@$ is a proper subquotient of $L$. Conversely, an $\Cf$-linear functor from $\Cf\CC_p^L$ to $\gMod{\Cf}$ which vanishes on $p$-groups $P$ such that $P^@\ncong L$ can be viewed as an $\Cf$-linear functor from $\Cf\CC_p^{\sharp L}$ to $\gMod{\Cf}$. In the sequel, we will freely identify those two types of functors, and consider $\mathsf{Fun}_\Cf\big(\Cf\CC_p^{\sharp L},\gMod{\Cf}\big)$ as the full subcategory of $\mathsf{Fun}_\Cf\big(\Cf\CC_p^L,\gMod{\Cf}\big)$ consisting of functors which vanish on $p$-groups $P$ such that $P^@\ncong L$.
\end{rem}
\begin{mth}{Theorem}\label{equivalencediese}{\rm [ }$p\in \Cf^\times${\rm ] } Let $L$ be an atoric $p$-group. 
\begin{enumerate}
\item If $F$ is a $p$-biset functor over $\Cf$ such that $F=\widehat{b}_LF$, and $P$ is a finite $p$-group such that $L\not\sqsubseteq P$, then $F(P)=\zero$.
\item If $G$ is an $\Cf$-linear functor from $\Cf\CC_p^{\sharp L}$ to $\gMod{\Cf}$, then $\widehat{b}_L\CR_{\CY_L}(G)=\CR_{\CY_L}(G)$.
\item The forgetful functor $\CO_{\CY_L}$ and its right adjoint $\CR_{\CY_L}$ restrict to quasi-inverse equivalences of categories
$$\xymatrix{
\widehat{b}_L\CF_{p,\Cf}\ar[r]<.5ex>^-{\CO_{\CY_L}}&\ar[l]<.5ex>^-{\CR_{\CY_L}} \mathsf{Fun}_\Cf\big(\Cf\CC_p^{\sharp L},\gMod{\Cf}\big)\mpoint
}$$
\end{enumerate}
\end{mth}
\pf (1) Since $\widehat{b}_LF=F$, then in particular $F(b_L^P)F(P)=F(P)$. If $L\not\sqsubseteq P$, then there is no minimal section $(T,S)$ of $P$ with $(T/S)^@\cong L$, thus $b_L^P=0$, and $F(P)=\zero$.\mpn
(2) Let $G$ be an $\Cf$-linear functor from $\Cf\CC_p^{\sharp L}$ to $\gMod{\Cf}$, in other words an $\Cf$-linear functor from $\CF\CC_p^L$ to $\gMod{\Cf}$ which vanishes on $p$-groups $P$ such that $P^@$ is a proper subquotient of $L$. By Theorem~\ref{equivalence}, we have $\widehat{b}_L^+\CR_{\CY_L}(G)=\CR_{\CY_L}(G)$. If $H$ is an atoric $p$-group which is a proper subquotient of $L$, then $G$ vanishes over any subquotient $Q$ of $H$, since $Q^@\sqsubseteq H\sqsubset L$ if $Q\sqsubseteq H$. In particular $b_H^P$ acts by 0 on $\CR_{\CY_L}(G)(P)$, for any finite $p$-group $P$: indeed $b_H^P$ is a linear combination of terms of the form $\Indinf_{X/M}^P\Defres_{X/M}^P$, where $(X,M)$ is a section of $P$ such that $S\leq M\leq\Phi(T)\leq X\leq T$, for some section $(T,S)$ of $P$ with $(T/S)^@\cong H$. For such a section $(X,M)$ of $P$, we have $(X/M)^@\sqsubseteq (T/S)^@\sqsubseteq H$, thus $G$ vanishes on any subquotient of $X/M$, so $\CR_{\CY_L}(G)(X/M)=\zero$, hence $b_H^P=0$ on $\CR_{\CY_L}(G)(P)$, as claimed. It follows that $\widehat{b}_H\CR_{\CY_L}(G)=0$, hence $\widehat{b}_L^+\CR_{\CY_L}(G)=\CR_{\CY_L}(G)=\widehat{b}_L\CR_{\CY_L}(G)$.\mpn
(3) This is a straightforward consequence of (1) and (2), by Theorem~\ref{equivalence}.\endpf
The following proposition gives some detail on the structure of the category $\Cf\CC_p^{\sharp L}$:
\begin{mth}{Proposition}\label{basis2} Let $p$ be a prime, and $L$ be an atoric $p$-group.
\begin{enumerate}
\item Let $P$ be a finite $p$-group. Then $P^@\cong L$ if and only if there exists an elementary abelian $p$-group $E$ such that $P\cong E\times L$.
\item Let $P=E\times L$ and $Q=F\times L$, where $E$ and $F$ are elementary abelian $p$-groups. If $M\leq (Q\times P)$, then $q(M)^@\cong L$ if and only if 
$$p_{1,2}(M)=p_{2,2}(M)=L\;\;\;\hbox{and}\;\;\;k_{1,2}(M)=k_{2,2}(M)=\un\mvirg$$
where $p_{1,2}$ and $p_{2,2}$ are the morphisms from $\big((H\times L)\times (G\times L)\big)$ to~$L$ defined by $p_{1,2}\big((h,x),(g,y)\big)=x$ and $p_{2,2}\big((h,x),(g,y)\big)=y$, and
\begin{eqnarray*}
k_{1,2}(M)&=&\{x\in L\mid \big((1,x),(1,1)\big)\in M\}\mvirg\\
k_{2,2}(M)&=&\{x\in L\mid \big((1,1),(1,x)\big)\in M\}\mpoint\\
\end{eqnarray*}
\end{enumerate}
\end{mth}
\pf (1) This follows from Proposition~\ref{ExL}.\mpn
(2) By Lemma~~\ref{abelemxL}, the $\Cf$-module $\Cf B(Q,P)$ has a basis consisting of the isomorphism classes of $(Q,P)$-bisets of the form $(Q\times P)/M$, where $M$ is a subgroup of $(Q\times P)$, up to conjugation, and $q(M)^@\cong L$. If $M$ is such a subgroup, then $L\cong \big(p_1(M)/k_1(M)\big)^@\sqsubseteq \big(p_1(M)\big)^@\sqsubseteq Q^@\cong L$, so $p_1(M)^@\cong\nolinebreak L$, and similarly $p_2(M)^@\cong L$. By Proposition~\ref{ExL} $p_1(M)^@\cong L$ if and only if $Ep_1(M)=\nolinebreak P$, which in turn is equivalent to $p_{1,2}(M)=L$. Similarly $p_2(M)^@\cong L$ if and only if $p_{2,2}(M)=L$. \par
Then $\big(p_1(M)/k_1(M)\big)^@\cong L$ if and only if $k_1(M)\cap\Phi\big(p_1(M)\big)=\un$, by Proposition~\ref{PatisoQat}. Moreover $\Phi\big(p_1(M)\big)=\Phi(P)$, as there exists an elementary abelian subgroup $E'$ of $P$ such that $P=E'\times p_1(M)$, by Proposition~\ref{ExL} again. Since $\Phi(P)=\un\times\Phi(L)$, it follows that $k_1(M)\cap \big(\un\times\Phi(L)\big)=\un$. Now $N=k_1(L)\cap (\un\times L)$ is a normal subgroup of $(\un\times L)$ (since $p_{1,2}(M)=\nolinebreak L$), which intersect trivially $\big(\un\times\Phi(L)\big)$. Since $L$ is atoric, by Lemma~\ref{atoric}, any central element of order~$p$ of $(\un\times L)$ is contained in $\big(\un\times\Phi(L)\big)$, so $N$ contains no non trivial central element of $(\un\times L)$, hence $N=\un$. Thus $k_1(L)\cap (\un\times L)=\un$, or equivalently $k_{1,2}(M)=\un$. Similarly $k_{2,2}(M)=\un$. Hence $q(M)^@\cong L$ if and only if $p_{1,2}(M)=p_{2,2}(M)=L$ and $k_{1,2}(M)=k_{2,2}(M)=\un$.\endpf
\section{$L$-enriched bisets}
\begin{mth}{Notation} Let $G$ and $H$ be finite groups. If $U$ is an $(H,G)$-biset, and $u\in U$, let $(H,G)_u$ denote the stabilizer of $u$ in $(H\times G)$, i.e.
$$(H,G)_u=\{(h,g)\in (H\times G)\mid hu=ug\}\mpoint$$
Let $H_u=k_1\big((H,G)_u\big)$ denote the stabilizer of $u$ in $H$, and $_uG=k_2\big((H,G)_u\big)$ denote the stabilizer of $u$ in $G$. Set moreover 
$$q(u)=q\big((H,G)_u\big)=(H,G)_u/(H_u\times {_uG})\mpoint$$
\end{mth}
\begin{mth}{Definition} Let $L$ be a finite group. For two finite groups $G$ and~$H$, an {\em $L$-enriched} $(H,G)$-biset is a $(H\times L,G\times L)$-biset $U$ such that $L\sqsubseteq q(u)$, for any $u\in U$. A morphism of $L$-enriched $(H,G)$-bisets is a morphism of $(H\times L,G\times L)$-bisets.\par
The disjoint union of two $L$-enriched $(H,G)$-bisets is again an $L$-enriched $(H,G)$-biset. Let $B\enr{L}(H,G)$ denote the Grothendieck group of finite $L$-enriched $(H,G)$-bisets for relations given by disjoint union decompositions. The group $B\enr{L}(H,G)$ is called the {\em Burnside group of $L$-enriched $(H,G)$-bisets}.  
\end{mth}
\begin{mth}{Lemma} Let $G, H, L$ be finite groups, and $U$ be an $(H\times L,G\times L)$-biset. Let $U^{\sharp L}$ denote the set of elements $u\in U$ such that $L\sqsubseteq q(u)$. Then $U^{\sharp L}$ is the largest sub-$L$-enriched $(H,G)$-biset of $U$.
\end{mth}
\pf It suffices to show that $U^{\sharp L}$ is a sub-$(H\times L,G\times L)$-biset of $U$, for then it is clearly the largest sub-$L$-enriched $(H,G)$-biset of $U$. And this is straightforward, since for any $(u,g,h,x,y)\in (U\times G\times H\times L\times L)$, if $v=(h,y)u(g,x)^{-1}$, then 
$$(H\times L,G\times L)_v={^{((h,y),(g,x))}(H\times L,G\times L)_u}\mvirg$$
and this conjugation induces a group isomorphism $q(v)\cong q(u)$.\endpf
\begin{mth}{Lemma}\label{basis} Let $G,H,L$ be finite groups.
\begin{enumerate}
\item Let $U$ be an $L$-enriched $(H,G)$-biset. If $V$ is a sub-$(H\times L,G\times L)$-biset of $U$, then $V$ is an $L$-enriched $(H,G)$-biset.
\item The group $B[L](H,G)$ has a $\Z$-basis consisting of the transitive bisets $\big((H\times L)\times (G\times L)\big)/M$, where $M$ is a subgroup of $\big((H\times L)\times (G\times L)\big)$ (up to conjugation) such that $L\sqsubseteq q(M)$.
\end{enumerate}
\end{mth}
\pf (1) This is straightforward.\mpn
(2) It follows from (1) that $B[L](H,G)$ has a basis consisting of the isomorphism classes of $L$-enriched $(H,G)$-bisets which are transitive $(H\times L,G\times L)$-bisets. These are of the form $U=\big((H\times L)\times (G\times L)\big)/M$, for some subgroup $M$ of $\big((H\times L)\times (G\times L)\big)$. Now if $u$ is the element $\big((1,1),(1,1)\big)M$ of $U$, the group $(H\times L,G\times L)_u$ is equal to $M$, hence $q(u)\cong q(M)$.\endpf
\begin{mth}{Lemma} Let $G,H,K,L$ be finite groups.
\begin{enumerate}
\item For an $(H,G)$-biset $U$, endow $U\times L$ with the $(H\times L,G\times L)$-biset structure defined by
$$\forall h\in H,\forall g\in G,\forall x,y,z\in L, \forall u\in U,\;\;\;(h,x)(u,y)(g,z)=(hug,xyz)\mpoint$$
Then $U\times L$ is an $L$-enriched $(H,G)$-biset.
\item In particular, for any finite group $G$, the identity biset of $G\times L$ is an $L$-enriched $(G,G)$-biset.
\item If $U$ is an $(H,G)$-biset and $V$ is a $(K,H)$-biset, then there is an isomorphism
$$(V\times L)\times_{(H\times L)}(U\times L)\cong (V\times_HU)\times L$$
of $L$-enriched $(H,G)$-bisets.
\end{enumerate}
\end{mth}  
\pf (1) For $u\in U$ and $l\in L$, 
$$(H\times L,G\times L)_{(u,l)}=\{\big((h,{^lx}),(g,x)\big)\mid hug=u,\;l\in L\}\cong (H,G)_u\times L\mpoint$$
In particular $(H\times L)_{(u,l)}=H_u\times\un$ and $_{(u,l)}(G\times L)={_uG}\times \un$, and $q\big((u,l)\big)\cong q(u)\times L$ has a (sub)quotient isomorphic to $L$.  \mpn
(2) In particular, if $H=G$ and $U$ is the identity $(G,G)$-biset, then $U\times L$ is the identity biset of $(G\times L)$.\mpn
(3) It is straightforward to check that the maps
$$\xymatrix@R=1ex{
[(v,x),(u,y)]\in (V\times L)\times_{(H\times L)}(U\times L)\ar@{|->}[r]&([v,u],xy)\in (V\times_H U)\times L\\
[(v,1),(u,l)]\in (V\times L)\times_{(H\times L)}(U\times L)&\ar@{|->}[l]([v,u],l)\in (V\times_H U)\times L
}
$$
are well defined isomorphisms of $(K\times L,G\times L)$-bisets, inverse to one another.\endpf
\newcommand{\timesup}[1]{{\mathop{\times}^{_{#1}}\limits}}
\begin{mth}{Notation}\label{timesup} Let $G,H,K,L$ be finite groups. If $U$ is an $L$-enriched $(H,G)$-biset and $V$ is an $L$-enriched $(K,H)$-biset, let $V\timesup{L}_HU$ denote the $L$-enriched $(K,G)$-biset defined by
$$V\timesup{L}_HU=\big(V\times_{(H\times L)}U\big)^{\sharp L}\mpoint$$ 
\end{mth}
\begin{mth}{Lemma} Let $G,H,J,K,L$ be finite groups.
\begin{enumerate}
\item If $V$ is a $(K\times L,H\times L)$-biset and $U$ is an $(H\times L,G\times L)$-biset, then
$$(V\times_{(H\times L)}U)^{\sharp L}=V^{\sharp L}\timesup{L}_HU^{\sharp L}\mpoint$$
In particular, if $V$ and $U$ are $L$-enriched bisets, so is $V\timesup{L}_HU$.
\item If $U$ and $U'$ are $L$-enriched $(H,G)$-bisets, if $V,V'$ are $L$-enriched $(K,H)$-bisets, then there are isomorphisms
\begin{eqnarray*}
V\timesup{L}_H(U\sqcup U')&\cong& (V\timesup{L}_HU)\sqcup(V\timesup{L}_HU')\\
(V\sqcup V')\timesup{L}_HU&\cong& (V\timesup{L}_HU)\sqcup(V'\timesup{L}_HU)
\end{eqnarray*}
of $L$-enriched $(K,G)$-bisets.
\item If moreover $W$ is an $L$-enriched $(J,K)$-biset, then there is a canonical isomorphism
$$(W\timesup{L}_KV)\timesup{L}_HU\cong W\timesup{L}_K(V\timesup{L}_HU)$$
of $L$-enriched $(J,G)$-bisets.
\end{enumerate}
\end{mth}
\pf (1) Denote by $[v,u]$ the image in $V\times_{(H\times L)}U$ of a pair $(v,u)\in (V\times U)$. By Lemma~2.3.20 of~\cite{bisetfunctors}, 
$$(K\times L,G\times L)_{[v,u]}=(K\times L,H\times L)_{v}*(H\times L,G\times L)_{u}\mvirg$$
so by Lemma~2.3.22 of~\cite{bisetfunctors}, the group $q\big([v,u]\big)$ is a subquotient of $q(v)$ and $q(u)$. So if $[v,u]\in (V\times_{(H\times L)}U)^{\sharp L}$, then $L$ is a subquotient of $q\big([v,u]\big)$, hence it is a subquotient of $q(v)$ and $q(u)$, that is $v\in V^{\sharp L}$ and $u\in U^{\sharp L}$. Hence
$$(V\times_{(H\times L)}U)^{\sharp L}\subseteq (V^{\sharp L}\times_{(H\times L)}U^{\sharp L})^{\sharp L}=V^{\sharp L}\timesup{L}_HU^{\sharp L}\mvirg$$
and the reverse inclusion $(V^{\sharp L}\times_{(H\times L)}U^{\sharp L})^{\sharp L}\subseteq (V\times_{(H\times L)}U)^{\sharp L}$ is obvious. Hence $(V\times_{(H\times L)}U)^{\sharp L}=V^{\sharp L}\timesup{L}_HU^{\sharp L}$. If $V$ and $U$ are $L$-enriched bisets, i.e. if $V=V^{\sharp L}$ and $U=U^{\sharp L}$, this gives $(V\times_{(H\times L)}U)^{\sharp L}=V\timesup{L}_HU$, so $V\timesup{L}_HU$ is an $L$-enriched biset.\mpn
(2) This is straightforward.\mpn
(3) With the above notation, there is a canonical isomorphism 
$$\alpha:(W\times_{(K\times L)}V)\times_{(H\times L)}U\to W\times_{(K\times L)}(V\times_{(H\times L)}U)$$
sending $\big[[w,v],u\big]$ to $\big[w,[v,u]\big]$. Hence
\begin{eqnarray*}
(W\timesup{L}_KV)\timesup{L}_HU&=&\big((W\timesup{L}_KV)\times_{(H\times L)}U\big)^{\sharp L}\\
&=&\big((W\times_{(K\times L)}V)^{\sharp L}\times_{(H\times L)}U\big)^{\sharp L}\\
&=&\big((W\times_{(K\times L)}V)\times_{(H\times L)}U\big)^{\sharp L}\;\;\;\hbox{[by (1)]}\\
\end{eqnarray*}
Similarly
\begin{eqnarray*}
W\timesup{L}_{K}(V\timesup{L}_HU)&=&\big(W\times_{(K\times L)}(V\timesup{L}_HU)\big)^{\sharp L}\\
&=& \big(W\times_{(K\times L)}(V\times_{(H\times L)}U)^{\sharp L}\big)^{\sharp L}\\
&=& \big(W\times_{(K\times L)}(V\times_{(H\times L)}U)\big)^{\sharp L}\;\;\;\hbox{[by (1)]}\mpoint\\
\end{eqnarray*}
Hence $\alpha$ induces an isomorphism $(W\timesup{L}_KV)\timesup{L}_HU\cong W\timesup{L}_K(V\timesup{L}_HU)$.  \endpf
\pagebreak[3]
\begin{mth}{Definition} Let $L$ be a finite group, and $\Cf$ be a commutative ring. The {\em $L$-enriched biset category} $\Cf\CC\enr{L}$ of finite groups over $\Cf$ is defined as follows:
\begin{itemize}
\item The objects of $\Cf\CC\enr{L}$ are the finite groups.
\item For finite groups $G$ and $H$,
$$\Hom_{\Cf\CC\enr{L}}(G,H)=\Cf\otimes_\Z B\enr{L}(H,G)=\Cf B\enr{L}(H,G)$$
is the $\Cf$-linear extension of the Burnside group of $L$-enriched $(H,G)$-bisets.
\item The composition in $\Cf\CC\enr{L}$ is the $\Cf$-linear extension of the product $(V,U)\mapsto V\timesup{L}_HU$ defined in~\ref{timesup}.
\item The identity morphism of the group $G$ is (image in $\Cf B\enr{L}(G,G)$ of) the identity biset of $G\times L$, viewed as an $L$-enriched $(G,G)$-biset.
\end{itemize}
The category $\Cf\CC\enr{L}$ is $\Cf$-linear. An {\em $L$-enriched biset functor} over $\Cf$ is an $\Cf$-linear functor from $\Cf\CC\enr{L}$ to $\gMod{\Cf}$. The category of $L$-enriched biset functors over $\Cf$ is denoted by ${\CF_{\Cf}}\enr{L}$. It is an abelian $\Cf$-linear category.
\end{mth}
\pagebreak[3]
\def\El{\CE{l}}
\begin{mth}{Theorem} Let $p$ be a prime number, and $\Cf$ be a commutative ring.
\begin{enumerate}
\item If $L$ is an atoric $p$-group, the category $\Cf\CC_p^{\sharp L}$ of Definition~\ref{CCpdeL} is equivalent to the full subcategory ${\Cf\El_p[L]}$ of $\Cf\CC\enr{L}$ consisting of elementary abelian $p$-groups.
\item If $p\in\CF^\times$, the category $\CF_{p,\Cf}$ of $p$-biset functors over $\Cf$ is equivalent to the direct product of the categories $\mathsf{Fun}_\Cf\big(\Cf\El_p[L],\gMod{\Cf}\big)$ of $\Cf$-linear functors from $\Cf\El_p[L]$ to $\gMod{\Cf}$, for $L\in[\At_p]$.
\end{enumerate}
\end{mth}
\pf (1) Let $E$ be an elementary abelian $p$-group. Then $(E\times L)^@\cong L$, so $E\times L$ is an object of $\Cf\CC_p^{\sharp L}$. Set $\CI(E)=E\times L$. If $E$ and $F$ are elementary abelian $p$-groups, and if $U$ is a finite $L$-enriched $(F,E)$-biset, then $U$ is in particular an $(F\times L,E\times L)$-biset, an we can consider its image $\CI(U)$ in the quotient $\Hom_{\Cf\CC_p^{\sharp L}}(E\times L,F\times L)$ of $\Cf B(F\times L,E\times L)$. This yields a unique $\Cf$-linear map $\Cf B[L](F,E)\to \Hom_{\Cf\CC_p^{\sharp L}}(E\times L,F\times L)$, still denoted by $\CI$. \par
We claim that these assignments define a functor $\CI$ from $\Cf\El_p[L]$ to $\Cf\CC_p^{\sharp L}$: indeed, the identity $(E\times L,E\times L)$-biset is clearly mapped to the identity morphism of $\CI(E)$. 
Moreover, if $G$ is an elementary abelian $p$-group, if $V$ is an $L$-enriched $(G,F)$-biset and $U$ is an $L$-enriched $(F,E)$-biset, it is clear that
$$\CI\big(V\timesup{L}_FU)=\CI(V)\circ \CI(U)\mvirg$$
where the right hand side composition is in the category $\Cf\CC_p^{\sharp L}$: indeed, the transitive bisets $(Q\times P)/M$ with $q(M)^@\sqsubset L$ appearing in the product $V\times_{(F\times L)}U$ are exactly those vanishing in $\Hom_{\Cf\CC_p^{\sharp L}}\big(\CI(E),\CI(F)\big)$, by Lemma~\ref{abelemxL}.
Hence $\CI$ is an isomorphism
$$\CI:\Cf B[L](F,E) \to \Hom_{\Cf\CC_p^{\sharp L}}\big(\CI(E),\CI(F)\big)\mpoint$$
In other words $\CI$ is a fully faithful functor from $\Cf\El_p[L]$ to $\Cf\CC_p^{\sharp L}$. Moreover, by Proposition~\ref{ExL}, if $P$ is a finite $p$-group with $P^@\cong L$, there exists an elementary abelian $p$-group $E$ such that $P$ is isomorphic to $E\times L$, hence $P$ is isomorphic to $E\times L$ in the category $\Cf\CC_p^{\sharp L}$.\par
It follows that the functor $\CI$ is fully faithful and essentially surjective, so it is an equivalence of categories.\mpn
(2) This is a straightforward consequence of (1), Assertion~5 of Corollary~\ref{decomposition}, and Assertion~3 of Theorem~\ref{equivalencediese}.\endpf
\section{The category $\widehat{b}_L\CF_{p,\Cf}$, for an atoric $p$-group $L$ ($p\in R^\times$)}
Let $L$ be a fixed atoric $p$-group. In this section, we give some detail on the structure of the category $\widehat{b}_L\CF_{p,\Cf}$ of $p$-biset functors invariant by the idempotent $\widehat{b}_L$. \par
We start by straightforward consequences of Theorem~\ref{equivalencediese}. For a finite $p$-group $P$, we denote by $\Sigma_{\sharp L}(P)$ the subset of $\Sigma_L(P)$ consisting of sections $(X,M)$ of $P$ such that $(X/M)^@\cong L$. When $G$ is an $\Cf$-linear functor from $\Cf\CC_p^{\sharp L}$ to $\gMod{\Cf}$, we can compute $\CR_{\CY_L}(G)$ at $P$ by restricting the inverse limit of~\ref{inverse limit} to the subset $\Sigma_{\sharp L}(P)$, i.e. by
$$\CR_{\CY_L}(G)(P)=\limproj{(X,M)\in\Sigma_{\sharp L}(P)}G(X/M)\mpoint$$
\begin{mth}{Proposition}\label{isofaciles} {\rm [ }$p\in \Cf^\times${\rm ] } Let $L$ be an atoric $p$-group. If $F$ is a $p$-biset functor in $\widehat{b}_L\CF_{p,\Cf}$, and $P$ is a finite $p$-group, then
\begin{eqnarray*}
F(P)&\cong& \limproj{(X,M)\in\Sigma_{\sharp L}(P)}F(X/M)\mvirg\\
&\cong&\dirsum{\substack{(T,S)\in[\CM(P)]\\(T/S)^@\cong L}}\delta_\Phi F(T/S)^{N_P(T,S)/T}\mpoint
\end{eqnarray*}
\end{mth}
\pf The isomorphism $F(P)\cong\limproj{(X,M)\in\Sigma_{\sharp L}(P)}F(X/M)$ is Assertion~3 of Theorem~\ref{equivalencediese}. The second isomorphism follows from Theorem~\ref{evaluation}, which implies that for $(T,S)\in \CM(P)$
$$\delta_\Phi F(T/S)^{N_P(T,S)/T}\cong F(\epsilon_{T,S}^P)\big(F(P)\big)\mpoint$$
Moreover $F(b_L^P)F(P)=F(P)$ since $F\in\widehat{b}_L\CF_{p,\Cf}$, and 
$$F(\epsilon_{T,S}^P)F(b_L^P)=F(\epsilon_{T,S}^Pb_L^P)=0$$
unless $(T/S)^@\cong L$. Thus $\delta_\Phi F(T/S)^{N_P(T,S)/T}=\zero$ unless $(T/S)^@\cong L$, which completes the proof.
\endpf
\pagebreak[3]
The decomposition of the category $\CF_{p,\Cf}$ of $p$-biset functors stated in Corollary~\ref{decomposition} leads to the following natural definition:
\begin{mth}{Definition}\label{vertex}{\rm [ }$p\in \Cf^\times${\rm ] }  Let $F$ be an indecomposable $p$-biset functor over~$\Cf$. There exists a unique atoric $p$-group $L$ (up to isomorphism) such that $F=\widehat{b}_LF$. The group $L$ is called the {\em vertex} of $F$.
\end{mth}
\begin{rem}{Remark} \label{vertex and ext}It follows in particular from this definition that if $F$ and~$F'$ are indecomposable $p$-biset functors over $\Cf$ with non-isomorphic vertices, then $\Ext_{\CF_{p,\Cf}}^*(F,F')=\zero$.
\end{rem}
\begin{mth}{Theorem}\label{vertex minimal}{\rm [ }$p\in \Cf^\times${\rm ] } Let $F$ be an indecomposable $p$-biset functor over~$\Cf$ and let $L$ be a vertex of $F$. If $Q$ is a finite $p$-group such that $F(Q)\neq \zero$, but $F$ vanishes on any proper subquotient of $Q$, then $L\cong Q^@$.
\end{mth}
\pf Let $Q$ be a finite $p$-group such that $F(Q)\neq\zero$ and $F(Q')=\zero$ for any proper subquotient $Q'$ of $Q$. By Proposition~\ref{formule epsilon}, if $(T,S)$ is a minimal section of $Q$, then
$$\epsilon_{T,S}^Q\!=\frac{1}{|N_Q(T,S)|}\!\!\!\!\!\!\sumb{X\leq T, M\normal T}{\rule{0ex}{1.6ex}S\leq M\leq\Phi(T)\leq X\leq T}|X|\mu(X,T)\mu_{\normal T}(S,M)\,\Indinf_{X/M}^Q\circ\Defres_{X/M}^Q\mpoint$$
Now if $X/M$ is a proper subquotient of $Q$, i.e. if $X\neq Q$ or $M\neq \un$, then $F(X/M)=\zero$, and $F(\Indinf_{X/M}^Q\circ\Defres_{X/M}^Q)=0$. Hence $F(\epsilon_{T,S}^Q)=0$ unless $T=Q$ and $S=\un$, and moreover
$$F(\epsilon_{Q,\un}^Q)=\frac{1}{|Q|}|Q|\mu(Q,Q)\mu_{\normal Q}(\un,Q)F(\Indinf_{Q/\un}^Q\Defres_{Q/\un}^Q)=\Id_{F(Q)}\mpoint$$
If $\widehat{b}_LF=F$, then in particular $F(b_L^Q)$ is equal to the identity map of $F(Q)$. This can only occur if the idempotent $\epsilon_{Q,\un}^Q$ appears in the sum defining $b_L^Q$, in other words if $(Q/\un)^@\cong L$, i.e. $Q^@\cong L$. Conversely, if $Q^@\cong L$, then $F(b_L^Q)=F(\epsilon_{Q,\un}^Q)=\Id_{F(Q)}\neq 0$. It follows that $\widehat{b}_LF\neq 0$, hence $\widehat{b}_LF=F$, since $F$ is indecomposable. Hence $Q^@$ is (isomorphic to) the vertex of $F$, as was to be shown.\endpf
We assume from now on that $\Cf=k$ is a field.\def\Cf{k} Recall (\cite{bisetfunctors} Chapter~4) that the simple $p$-biset functors over $\Cf$ are indexed by pairs $(Q,V)$ consisting of a $p$-group $Q$ and a simple $\Cf\Out(Q)$-module $V$.
\pagebreak[3]
\begin{mth}{Corollary}\label{vertex simple} Let $\Cf$ be a field of characteristic different from $p$. 
\begin{enumerate}
\item If $Q$ is a finite $p$-group, and $V$ is a simple $\Cf \Out(Q)$-module, then the vertex of the simple $p$-biset functor $S_{Q,V}$ is isomorphic to $Q^@$.
\item Let $Q$ (resp. $Q'$) be a finite $p$-group, and $V$ (resp. $V'$) be a simple $k\Out(Q)$-module (resp. a simple $k\Out(Q')$-module). If $Q^@\ncong Q'^@$, then $\Ext_{\CF_{p,k}}^*(S_{Q,V},S_{Q',V'})=\zero$.
\end{enumerate}
\end{mth}
\pf (1) Indeed $Q$ is a minimal group for $S_{Q,V}$, so $S_{Q,V}(Q)\neq \zero$, but $S_{Q,V}$ vanishes on any proper subquotient of $Q$.\mpn
(2) Follows from (1) and Remark~\ref{vertex and ext}.\endpf
\pagebreak[3]
\begin{mth}{Definition} Let $F$ be a $p$-biset functor. A functor $S$ is a {\em subquotient} of $F$ (notation $S\sqsubseteq F$) if there exist subfunctors $F_2<F_1\leq F$ such that $F_1/F_2\cong S$. A {\em composition factor} of $F$ is a simple subquotient of $F$.\end{mth}
\begin{mth}{Lemma}\label{greatest} Let $\Cf$ be a field, and $F$ be a $p$-biset functor over $\Cf$.
\begin{enumerate}
\item If $F$ is a non zero, then $F$ admits a composition factor. 
\item If $\mathcal{S}$ is a family of simple $p$-biset functors over $\Cf$, there exists a greatest subfunctor of $F$ all composition factors of which belong to~$\mathcal{S}$.
\end{enumerate}
\end{mth}
\pf (1) Let $P$ be a finite $p$-group such that $F(P)\neq\zero$. Then $F(P)$ is a $\Cf B(P,P)$-module. Choose $m\in F(P)-\zero$, and consider the $\Cf B(P,P)$-submodule $M$ of $F(P)$ generated by $m$. Since $\Cf B(P,P)$ is finite dimensional over~$\Cf$, the module $M$ is also finite dimensional over $\Cf$, hence it contains a simple submodule $V$. By Proposition~3.1 of~\cite{boustathe}, there exists a simple $p$-biset functor $S$ such that $S(P)\cong V$ as $\Cf B(P,P)$-module. Then $S(P)$ is a subquotient of $F(P)$, so by Proposition~3.5 of~\cite{boustathe}, there exists a subquotient of $F$ isomorphic to $S$.\mpn
(2) Observe first that if $M,N$ are subfunctors of $F$, then any composition factor of $M+N$ is a composition factor of $M$ or a composition factor of $N$: indeed, if $S$ is a composition factor of $M+N$, let $F_2<F_1\leq M+N$ with $S\cong F_2/F_1$, and consider the images $F'_1$ and $F'_2$ of $F_1$ and $F_2$, respectively,  in the quotient $(M+N)/N\cong M/(M\cap N)$. If $F'_1\neq F'_2$, that is if $F_1+N\neq F_2+N$, then $F'_1/F'_2\cong (F_1+N)/(F_2+N)\cong F_1/F_2\cong S$ is a subquotient of $(M+N)/N\cong M/(M\cap N)$, hence $S$ is a subquotient of $M$. Otherwise $F_1+N=F_2+N$, so $F_1=F_2+(F_1\cap N)$, hence $S\cong F_1/F_2\cong (F_1\cap N)/(F_2\cap N)$ is a subquotient of~$N$. It follows by induction that any subquotient $S$ of a finite sum $\sum_{M\in\CI}\limits M$ of subfunctors of $F$ is a subquotient of some $M\in\CI$.\par

The latter also holds when $\CI$ is infinite: let $\Sigma=\sum_{M\in \CI}\limits M$ be an arbitrary sum of subfunctors of $F$, and $S$ be a composition factor of $\Sigma$. Let $F_2<F_1$ be subfunctors of $\Sigma$ such that $S\cong F_1/F_2$. If $P$ is a $p$-group such that $S(P)\cong F_1(P)/F_2(P)\neq 0$, let $U$ be a finite subset of $F_1(P)$ such that $F_1(P)/F_2(P)$ is generated as a $\Cf B(P,P)$-module by the images of the elements of $U$ (such a set exists because $S(P)$ is finite dimensional over $\Cf$, for any $P$). If $V$ is the $\Cf B(P,P)$-submodule of $F_1(P)$ generated by $U$, then $V$ maps surjectively on the module $F_1(P)/F_2(P)$, so there is a $\Cf B(P,P)$-submodule $W$ of $V$ such that $V/W\cong S(P)$. Now since $U$ is finite, there exists a finite subset $\mathcal{J}$ of $\CI$ such that $U\subseteq \sum_{M\in\mathcal{J}}\limits M(P)$. Setting $\Sigma_1=\sum_{M\in\mathcal{J}}\limits M$, it follows that $V/W\cong S(P)$ is a subquotient of $\Sigma_1(P)$, so by Proposition~3.5 of~\cite{boustathe}, there exists a subquotient of $\Sigma_1$ isomorphic to $S$. By the observation above $S$ is a subquotient of some $M\in\CJ\subseteq \CI$.\par
Now let $\mathcal{I}$ the set of subfunctors $M$ of $F$ such that all the composition factors of $M$ belong to $\CS$, and $N=\sum_{M\in\CI}\limits M$. The above discussion shows that $N\in\CI$, so $N$ is the greatest element of $\CI$. 
\endpf
\begin{mth}{Theorem} Let $\Cf$ be a field of characteristic different from $p$, and $L$ be an atoric $p$-group. Let $\CF_{p,\Cf}[L]$ the full subcategory of $\CF_{p,\Cf}$ consisting of functors  whose composition factors all have vertex $L$, i.e. are all isomorphic to $S_{P,V}$, for some $p$-group $P$ such that $P^@\cong L$, and some simple $\Cf\Out(P)$-module~$V$. 
\begin{enumerate}
\item If $F$ is a $p$-biset functor, then $\widehat{b}_L\CF_{p,\Cf}$ is the greatest subfunctor of $F$ which belongs to $\CF_{p,\Cf}[L]$.
\item In particular $\widehat{b}_L\CF_{p,\Cf}=\CF_{p,\Cf}[L]$.
\end{enumerate}
\end{mth}
\pf (1) Let $F$ be a $p$-biset functor over $\Cf$, and let $F_1=\widehat{b}_LF$. If $S$ is a composition factor of $F_1$, then $S=\widehat{b}_LS$, as $S$ is a subquotient of $F_1$. Hence $S$ has vertex~$L$, by Definition~\ref{vertex}. It follows that $F_1$ is contained in the greatest subfunctor $F_2$ of $F$ which belongs to $\CF_{p,\Cf}[L]$ (such a subfunctor exists by Lemma~\ref{greatest}). \par
Conversely, we know that $F_2=\dirsum{Q\in[\At_p]}\widehat{b}_QF_2$. For $Q\in[\At_p]$, any composition factor $S$ of $\widehat{b}_QF_2$ has vertex $Q$, by Definition~\ref{vertex}. But $S$ is also a direct summand of $F_2$, so $Q\cong L$. It follows that if $Q\ncong L$, then $\widehat{b}_QF_2$ has no composition factor, so $\widehat{b}_QF_2=\zero$, by Lemma~\ref{greatest}. In other words $F_2=\widehat{b}_LF_2$, hence $F_2\leq F_1$, and $F_2=F_1$, as was to be shown.\mpn
(2) Let $F$ be a $p$-biset functor. Then $F\in \widehat{b}_L\CF_{p,\Cf}$ if and only if $F=\widehat{b}_LF$, i.e. by (1) if and only if all the composition factors of $F$ have vertex $L$.\endpf
\begin{rem}{Example} {\bf the Burnside functor.} Let $\Cf$ be a field of characteristic $q\neq p$ ($q\geq 0$). It was shown in~\cite{both} Theorem 8.2 (see also~\cite{bisetfunctors} 5.6.9) that the Burnside functor $\Cf B$ is uniserial, hence indecomposable. As $\Cf B(\un)\neq 0$, the vertex of~$\Cf B$ is the trivial group, by Theorem~\ref{vertex minimal}, thus $\Cf B$ is an object of $\widehat{b}_\un\CF_{p,\Cf}=\CF_{p,\Cf}[\un]$. It means that all the composition factors of $\Cf B$ have to be of form $S_{Q,V}$, where $Q^@=\un$, i.e. $Q$ is elementary abelian. And indeed by~\cite{both} Theorem 8.2, the composition factors of $\Cf B$ are all of the form $S_{Q,\Cf}$, where $Q$ runs through a specific set of elementary abelian $p$-groups which depends on the order of $p$ modulo $q$ (suitably interpreted when $q=0$).
\end{rem}

\centerline{\rule{2cm}{.1ex}}\vspace{2ex}
Serge Bouc\\
LAMFA-CNRS UMR7352\\
33 rue St Leu, 80039 - Amiens - Cedex 01\\
France\\
{\tt email: serge.bouc@u-picardie.fr}
\end{document}